\documentclass[a4paper]{article}
\pdfoutput=1

\usepackage{graphicx,wrapfig}
\usepackage{lipsum}
\usepackage{amsmath,amssymb,mathrsfs}
\usepackage{hyperref}
\usepackage{url}
\usepackage{mathtools}
\usepackage{changepage}
\usepackage{multirow}
\usepackage{wrapfig}
\usepackage{subfig}
\usepackage{color}
\usepackage{epsfig,enumerate,amsmath,amsfonts,amssymb,amsthm,mathrsfs,ifpdf}
\usepackage{indentfirst,relsize}
\usepackage{setspace,graphicx}
\usepackage{latexsym}
\usepackage[all]{xy}
\usepackage[usenames,dvipsnames]{pstricks}
\usepackage{pst-grad} % For gradients
\usepackage{pst-plot} % For axes
\usepackage[margin = 2.50cm]{geometry}

\usepackage{algorithm}

\usepackage{algpseudocode}

\newcommand{\floor}[1]{\left\lfloor #1 \right\rfloor}
\newcommand{\ceil}[1]{\left\lceil #1 \right\rceil}

\newcommand{\remove}[1]{}

\newtheorem{theorem}{Theorem}%[section]
\newtheorem{lemma}[theorem]{Lemma}
\newtheorem{proposition}[theorem]{Proposition}
\newtheorem{corollary}[theorem]{Corollary}

\newtheorem{claim}[theorem]{Claim}
\newcounter{Ca}[theorem]
\newtheorem{case}[Ca]{Case}
\newtheorem{definition}[theorem]{Definition}
\newtheorem{remark}{Remark}

\title{Algorithmic study on $2$-transitivity of graphs}
\author{Subhabrata Paul \and Kamal Santra}

\author{Subhabrata Paul\footnote{Department of Mathematics, IIT Patna, India, email:subhabrata@iitp.ac.in} \and Kamal Santra\footnote{Department of Mathematics, IIT Patna, India, email:kamal\_1821ma04@iitp.ac.in} }

\date{}

\begin{document}

\maketitle
\begin{abstract}
		Let $G=(V, E)$ be a graph where $V$ and $E$ are the vertex and edge sets, respectively. For two disjoint subsets $A$ and $B$ of $V$, we say $A$ \emph{dominates} $B$ if every vertex of $B$ is adjacent to at least one vertex of $A$. A vertex partition $\pi = \{V_1, V_2, \ldots, V_k\}$ of $G$ is called a \emph{transitive partition} of size $k$ if $V_i$ dominates $V_j$ for all $1\leq i<j\leq k$. In this article, we study a variation of transitive partition, namely \emph{$2$-transitive partition}. For two disjoint subsets $A$ and $B$ of $V$, we say $A$ \emph{$2$-dominates} $B$ if every vertex of $B$ is adjacent to at least two vertices of $A$. A vertex partition $\pi = \{V_1, V_2, \ldots, V_k\}$ of $G$ is called a \emph{$2$-transitive partition} of size $k$ if $V_i$ $2$-dominates $V_j$ for all $1\leq i<j\leq k$. The \textsc{Maximum $2$-Transitivity Problem} is to find a $2$-transitive partition of a given graph with the maximum number of parts. We show that the decision version of this problem is NP-complete for chordal and bipartite graphs. On the positive side, we design three linear-time algorithms for solving \textsc{Maximum $2$-Transitivity Problem} in trees, split and bipartite chain graphs.

\end{abstract}


	\section{Introduction}

Partitioning a graph is one of the fundamental problems in graph theory. In the partitioning problem, the objective is to partition the vertex set (or edge set) into some parts with desired properties, such as independence, minimal edges across partite sets, etc. In literature, partitioning the vertex set into certain parts so that the partite sets follow particular domination relations among themselves has been studied. Let $G$ be a graph with $V(G)$ as its vertex set and $E(G)$ as its edge set. When the context is clear, $V$ and $E$ are used instead of $V(G)$ and $E(G)$. The \emph{neighbourhood} of a vertex $v\in V$ in a graph $G=(V, E)$ is the set of all adjacent vertices of $v$ and is denoted by $N_G(v)$. The \emph{degree} of a vertex $v$ in $G$, denoted as $\deg_G(v)$, is the number of edges incident to $v$. A vertex $v$ is said to \emph{dominate} itself and all its neighbouring vertices. A \emph{dominating set} of $G=(V, E)$ is a subset of vertices $D$ such that every vertex $x\in V\setminus D$ has a neighbour $y\in D$, that is, $x$ is dominated by some vertex $y$ of $D$. For two disjoint subsets $A$ and $B$ of $V$, we say $A$ \emph{dominates} $B$ if every vertex of $B$ is adjacent to at least one vertex of $A$.

There has been a lot of research on graph partitioning problems based on a domination relationship between the different sets. Cockayne and Hedetniemi introduced the concept of \emph{domatic partition} of a graph $G=(V, E)$ in 1977, in which the vertex set is partitioned into $k$ parts, say $\pi =\{V_1, V_2, \ldots, V_k\}$, such that each $V_i$ is a dominating set of $G$ \cite{cockayne1977towards}. The number representing the highest possible order of a domatic partition is called the \emph{domatic number} of G, denoted by $d(G)$. Another similar type of partitioning problem is the \emph{Grundy partition}. Christen and Selkow introduced a Grundy partition of a graph $G=(V, E)$ in 1979 \cite{CHRISTEN197949}. In the Grundy partitioning problem, the vertex set is partitioned into $k$ parts, say $\pi =\{V_1, V_2, \ldots, V_k\}$, such that each $V_i$ is an independent set and for all $1\leq i< j\leq k$, $V_i$ dominates $V_j$. The maximum order of such a partition is called the \emph{Grundy number} of $G$, denoted by $\Gamma(G)$. In 2018, J. T. Hedetniemi and S. T. Hedetniemi \cite{hedetniemi2018transitivity} introduced a transitive partition as a generalization of the Grundy partition. A \emph{transitive partition} of size $k$ is defined as a partition of the vertex set into $k$ parts, say $\pi =\{V_1,V_2, \ldots, V_k\}$, such that for all $1\leq i< j\leq k$, $V_i$ dominates $V_j$. The maximum order of such a transitive partition is called the \emph{transitivity} of $G$ and is denoted by $Tr(G)$. Recently, in 2020, Haynes et al. generalized the idea of domatic partition as well as transitive partition and introduced the concept of \emph{upper domatic partition} of a graph $G$, where the vertex set is partitioned into $k$ parts, say $\pi =\{V_1, V_2, \ldots, V_k\}$, such that for each $i, j$, with $1\leq i<j\leq k$, either $V_i$ dominates $V_j$ or $V_j$ dominates $V_i$ or both \cite{haynes2020upper}. The maximum order of such an upper domatic partition is called the \emph{upper domatic number} of $G$, denoted by $D(G)$. All these problems, domatic number \cite{chang1994domatic, zelinka1980domatically, zelinka1983k}, Grundy number \cite{effantin2017note, furedi2008inequalities, hedetniemi1982linear, zaker2005grundy, zaker2006results}, transitivity \cite{haynes2019transitivity, hedetniemi2018transitivity, paul2023transitivity, santra2023transitivity}, upper domatic number \cite{haynes2020upper, samuel2020new} have been extensively studied both from an algorithmic and structural point of view. A Grundy partition is a transitive partition with the additional restriction that each partite set must be independent. In a transitive partition $\pi =\{V_1, V_2, \ldots, V_k\}$ of $G$, we require domination property in one direction, that is, $V_i$ dominates $V_j$ for $1\leq i< j\leq k$. However, in a upper domatic partition $\pi =\{V_1,V_2, \ldots, V_k\}$ of $G$, for all $1\leq i<j\leq k$, either $V_i$ dominates $V_j$ or $V_j$ dominates $V_i$ or both. The definition of each vertex partitioning problem ensures the following inequalities for any graph $G$. For any graph $G$, $1\leq \Gamma(G)\leq Tr(G)\leq D(G)\leq n$.

In this article, we introduce a similar graph partitioning problem, namely \emph{$2$-transitive partition}, which is a variation of transitive partition. For two disjoint subsets $A$ and $B$, we say $A$ \emph{$2$-dominates} $B$ if every vertex of $B$ is adjacent to at least two vertices of $A$. A \emph{$2$-transitive partition} of size $k$ is defined as a partition of the vertex set into $k$ parts, say $\pi =\{V_1,V_2, \ldots, V_k\}$, such that for all $1\leq i< j\leq k$, $V_i$ $2$-dominates $V_j$. The maximum order of such a $2$-transitive partition is called the \emph{$2$-transitivity} of $G$ and is denoted by $Tr_2(G)$. Note that every $2$-transitive partition is also a transitive partition. Therefore, for any graph $G$, $1\leq Tr_2(G)\leq Tr(G)\leq n$. Some vertex partition parameters exist where the parameter's value in a subgraph can be greater than the original graph. The upper domatic number is one such example. But in the case of a $2$-transitive partition, $Tr_{2}(H)\leq Tr_{2}(G)$, for every subgraph $H$ of $G$. Consequently, for a disconnected graph, the $2$-transitivity equals the maximum $2$-transitivity among its components. Therefore, we focus only on connected graphs in this paper. The \textsc{Maximum $2$-Transitivity Problem} and its corresponding decision version are defined as follows:

\noindent\textsc{\underline{Maximum $2$-Transitivity Problem(M$2$TP)}}

\noindent\emph{Instance:} A graph $G=(V,E)$

\noindent\emph{Solution:} An $2$-transitive partition of $G$

\noindent\emph{Measure:} Order of the $2$ transitive partition of $G$\\

\noindent\textsc{\underline{Maximum $2$-Transitivity Decision Problem(M$2$TDP)}}

\noindent\emph{Instance:} A graph $G=(V,E)$, integer $k$

\noindent\emph{Question:} Does $G$ have a $2$-transitive partition of order at least $k$?\\

In this paper, we study the computational complexity of the $2$-transitivity problem. The main contributions are summarized below:
\begin{enumerate}
	
	\item[1.] We show that the \textsc{M$2$TDP} is NP-complete for chordal and bipartite graphs.
	
	\item [2.] We show that the \textsc{M$2$TP} can be solved in linear time for trees, split graphs, and bipartite chain graphs.
\end{enumerate}

The rest of the paper is organized as follows. Section 2 contains basic definitions and notations that are followed throughout the article. This section also discusses the properties of $2$-transitivity of graphs. Section 3 shows that the \textsc{M$2$TDP} is NP-complete in chordal and bipartite graphs. In Section 4, we design three linear-time algorithms for solving \textsc{M$2$TP} in trees, split graphs, and bipartite chain graphs. Finally, Section 5 concludes the article.

\section{Preliminaries}

\subsection{Definitions and notations}

Let $G=(V, E)$ be a graph with $V$ and $E$ as its vertex and edge sets, respectively. A graph $H=(V', E')$ is said to be a \emph{subgraph} of a graph $G=(V, E)$ if and only if $V'\subseteq V$ and $E'\subseteq E$. For a subset $S\subseteq V$, the \emph{induced subgraph} on $S$ of $G$ is defined as the subgraph of $G$ whose vertex set is $S$ and edge set consists of all of the edges in $E$ that have both endpoints in $S$, and it is denoted by $G[S]$. The \emph{complement} of a graph $G=(V,E)$ is the graph $\overline{G}=(\overline{V}, \overline{E})$, such that $\overline{V}=V$ and $\overline{E}=\{uv| uv\notin E \text{ and } u\neq v\}$. The \emph{open neighbourhood} of a vertex $x\in V$ is the set of vertices $y$ adjacent to $x$, denoted by $N_G(x)$. The \emph{closed neighborhood} of a vertex $x\in V$, denoted as $N_G[x]$, is defined by $N_G[x]=N_G(x)\cup \{x\}$. 

A subset of $S\subseteq V$ is said to be an \emph{independent set} of $G$ if no two vertices in $S$ are adjacent. A subset of $K\subseteq V$ is said to be a \emph{clique} of $G$ if every pair of vertices in $K$ are adjacent. The cardinality of a maximum size clique is called \emph{clique number} of $G$, denoted by $\omega(G)$. 

A graph is called \emph{bipartite} if its vertex set can be partitioned into two independent sets. A bipartite graph $G=(X\cup Y,E)$ is called a \textit{bipartite chain graph} if there exists an ordering of vertices of $X$ and $Y$, say $\sigma_X= (x_1,x_2, \ldots ,x_{n_1})$ and $\sigma_Y=(y_1,y_2, \ldots ,y_{n_2})$, such that $N(x_{n_1})\subseteq N(x_{n_1-1})\subseteq \ldots \subseteq N(x_2)\subseteq N(x_1)$ and $N(y_{n_2})\subseteq N(y_{n_2-1})\subseteq \ldots \subseteq N(y_2)\subseteq N(y_1)$. The ordering of $X$ and $Y$ is called a \emph{chain ordering}, which can be computed in linear time \cite{heggernes2007linear}. A graph $G=(V, E)$ is a \emph{split graph} if $V$ can be partitioned into an independent set $S$ and a clique $K$.

An edge between two non-consecutive vertices of a cycle is called a \emph{chord}. If every cycle in $G$ of length at least four has a chord, then $G$ is called a \emph{chordal graph}. A vertex $v\in V$ is called a \emph{simplicial vertex} of $G$ if $N_G[v]$ induces a clique in $G$. A \emph{perfect elimination ordering (PEO)} of $G$ is an ordering of the vertices, say $\sigma=(v_1, v_2,\ldots, v_n)$, such that $v_i$ is a simplicial vertex of $G_i=G[\{v_i, v_{i+1}, \ldots , v_n\}]$ for all $1\leq i\leq n$. Chordal graphs can be characterized by the existence of PEO; that is, a graph $G$ is chordal if and only if $G$ has a PEO \cite{fulkerson1965incidence}.

\subsection{Properties of $2$-transitivity}
In this subsection, we discuss some properties of the $2$-transitivity of a graph $G$. We start by showing an upper bound for $2$-transitivity.

\begin{proposition} \label{upper_bound_2-transitivity}
	For any graph $G$, $Tr_2(G)\leq \floor{\frac{\Delta(G)}{2}}+1$, where $\Delta(G)$ is the maximum degree of $G$.
	
\end{proposition}

\begin{proof}
	Let $\pi= \{V_1, V_2, \ldots, V_k\}$  be a $Tr_2$-partition of $G$ and $x\in V_k$ be a vertex of $G$. Since $V_i$ $2$-dominates $V_k$ for all $1\leq i\leq k-1$, the degree of $x$ is at least $2(k-1)$. This implies that the maximum degree of $G$, that is, $\Delta(G)$ is at least $2(k-1)$. Therefore, we have $k\leq \floor{\frac{\Delta(G)}{2}}+1$ as $k$ is an integer.
\end{proof}

It can be easily verified that paths, cycles, and complete graphs on $n$ vertices achieve this bound, which is tight.
\begin{corollary}
	The following are true:
	\begin{itemize}
		\item[(a)] For a path $P_n$ with $n\geq 3$ vertices, $Tr_2(P_n)=2$.
		\item[(b)] For a cycle $C_n$ with $n\geq 3$ vertices, $Tr_2(C_n)=2$.
		\item[(c)] For a complete graph $K_n$ with $n$ vertices, $Tr_2(K_n)=\floor{\frac{n-1}{2}}+1$.
	\end{itemize}
\end{corollary}

%The $2$-transitivity of paths and cycles in the following proposition is immediately due to the above bound.
%
%\begin{pre}\label{Path_Tr_2}
%For the path $P_n$ and cycle $C_n$ with $n\geq 3$, $Tr_2(P_n)=Tr_2(C_n)=2$.
%\end{pre}
%
%
%For the complete graph $K_n$ an the complete bipartite graph $K_{m, n}$, we have the following propositions.
%
%\begin{pre}\label{complete_graph_2-transitivity}
%	If $G$ is a $K_n$, then $Tr_2(G)=\floor{\frac{n-1}{2}}+1$. 
%	
%\end{pre}
%
%\begin{proof}
%Let $n$ be odd and consider $\pi =\{V_1, V_2, \ldots, V_k\}$ vertex partition of $G$ such that each $V_i$, $1\leq i\leq k-1 $ contains exactly two vertices from $K_n$ and $V_k$ contains only one vertex from $K_n$. Clearly $\pi$ is a $2$-transitive partition of $G$ of order $k=\frac{n-1}{2}+1$. Therefore, $Tr_2(K_n)\geq \frac{n-1}{2}+1=\floor{\frac{n-1}{2}}+1$. But we know by the Proposition \ref{upper_bound_2-transitivity}, $Tr_2(K_n)\leq \floor{\frac{\Delta(G)}{2}}+1= \floor{\frac{n-1}{2}}+1$. Now, let $n$ be even and consider $\pi =\{V_1, V_2, \ldots, V_k\}$ vertex partition of $G$ such that each $V_i$, $1\leq i\leq k$ contains exactly two vertices from $K_n$. Clearly $\pi$ is a $2$-tarnsitive partition of $G$ of order $k=\frac{n}{2}$. Therefore, $Tr_2(K_n)\geq \frac{n}{2}=\floor{\frac{n-1}{2}}+1$. Similarly, $Tr_2(K_n)\leq \floor{\frac{\Delta(G)}{2}}+1= \floor{\frac{n-1}{2}}+1$. Hence, If $G$ is a $K_n$, then $Tr_2(G)=\floor{\frac{n-1}{2}}+1$.
%	
%\end{proof}

Next, we show a lower bound on $2$-transitivity of a graph $G$ in terms of the transitivity of $G$.

\begin{proposition} \label{relation_between_2-transitive_and_transitive}
	Let $G$ be a graph with $Tr(G)=k$. Then the $2$-transitivity of a graph $G$ is bounded as follows:
	\begin{itemize}
		\item[(a)] $Tr_2(G)\geq \frac{k}{2}$ if $k$ is even and
		\item[(b)] $Tr_2(G)\geq \frac{k+1}{2}$ if $k$ is odd.
	\end{itemize}
\end{proposition}

\begin{proof}
	Let $\pi =\{V_1, V_2, \ldots, V_k\}$ be a transitive partition of $G$ of size $k$. If $k$ is even, then 
	let $\pi'=\{U_1, U_2, \ldots, U_{\frac{k}{2}}\}$ be a vertex partition of $G$, such that $U_i= V_{2i-1}\cup V_{2i}$ for all $1\leq i \leq  \frac{k}{2}$. Consider two sets $U_i$ and $U_j$ such that $i<j$. Let $x\in U_j=V_{2j-1}\cup V_{2j} $. As $\pi$ is a transitive partition of $G$, there exists $y_1\in V_{2i-1}$ and $y_2\in V_{2i}$ such that $xy_1, xy_1\in E(G)$. So, $U_i$ $2$-dominates $U_j$. Therefore, $\pi'$ is a $2$-transitive partition and hence, $Tr_2(G)\geq \frac{k}{2}$. If $k$ is odd, then let $\pi'=\{U_1, U_2, \ldots, U_{\frac{k-1}{2}}, U_{\frac{k+1}{2}}\}$ be a vertex partition of $G$, such that $U_i= V_{2i-1}\cup V_{2i}$ for all $1\leq i \leq  \frac{k-1}{2}$ and $U_{\frac{k+1}{2}}=V_k$. Using similar arguments, we can show that $\pi'$ is a $2$-transitive partition and hence, $Tr_2(G)\geq \frac{k+1}{2}$. 
\end{proof}

The following result shows that the lower bound is achieved for certain graphs; therefore, the above bound is tight.

\begin{proposition}\label{2-transitivity_equal_transitivity/2}
	Let $G$ be a graph with $Tr(G)=\Delta(G)+1$, where $\Delta(G)$ is the maximum degree of the graph. Then 
	\begin{itemize}
		\item[(a)] $Tr_2(G)=\frac{Tr(G)}{2}$, when $Tr(G)$ is even and
		\item[(b)] $Tr_2(G) =\frac{Tr(G)+1}{2}$, when $Tr(G)$ is odd.
	\end{itemize}
	
\end{proposition}
\begin{proof}
	Let us assume $Tr(G)$ is even. So, $\Delta(G)$ is odd. Now, by Proposition \ref{upper_bound_2-transitivity}, we have $Tr_2(G)\leq \floor{\frac{\Delta(G)}{2}}+1$. So, we can write $Tr_2(G)\leq \frac{\Delta(G)-1}{2}+1$, as $\Delta(G)$ is odd. Therefore, $Tr_2(G)\leq \frac{Tr(G)-1-1}{2}+1=\frac{Tr(G)}{2}$. Also, by Proposition \ref{relation_between_2-transitive_and_transitive} we know that $Tr_2(G)\geq \frac{Tr(G)}{2}$. Hence, $Tr_2(G)=\frac{Tr(G)}{2}$. Similarly, when $Tr(G)$ is odd, by Proposition \ref{upper_bound_2-transitivity}, we have $Tr_2(G)\leq \floor{\frac{\Delta(G)}{2}}+1=\frac{\Delta(G)}{2}+1=\frac{Tr(G)-1}{2}+1=\frac{Tr(G)+1}{2}$ and by Proposition \ref{relation_between_2-transitive_and_transitive}, we have $Tr_2(G)\geq \frac{Tr(G)+1}{2}$. Hence, $Tr_2(G)=\frac{Tr(G)+1}{2}$.
\end{proof}

%The following Figure \ref{fig:converse_not trure} describe this things.
%
%\begin{figure}[htbp!]
%	\centering
%	\includegraphics[scale=0.75]{split_graph}
%	\caption{Split graph $G$, $I=\{v_1, v_2, v_3\}$ and $C=\{u_1, u_2, u_3, u_4\}$. $Tr_2(G)=2=\frac{Tr(G)}{2}$, $\Delta(G)+1=4+1\not=4=Tr(G)$.}
%	\label{fig:converse_not trure}
%\end{figure}

Next, to show that the difference, $Tr_2(G)- \frac{Tr(G)}{2}$ or $Tr_2(G)- \frac{Tr(G)+1}{2} $ can be arbitrarily large. We show that there are graphs with equal transitivity and $2$-transitivity. To this end, we define a special tree, namely \emph{$2$-complete minimum broadcast tree of order $k$}. The recursive definition is as follows:

\begin{definition}
	A $2$-complete minimum broadcast tree of order $k$ $($$2$-$cmbt_k$$)$ is defined recursively as follows:
	\begin{enumerate}
		\item The $2$-cmbt of order $1$ $($$2$-$cmbt_1$$)$ is a rooted tree with one vertex.
		
		%\item The $2$-cmbt of order $2$ $($$2$-$cmbt_2$$)$ consists of a rooted tree, whose root has degree $2$ and is adjacent to two $1$ orders of $2$-cmbt.
		
		\item The $2$-cmbt of order $k$ $($$2$-$cmbt_k$$)$ is a rooted tree, whose root has degree $2(k-1)$ and is adjacent to vertices $\{v_1, v_2, \ldots, v_{2(k-1)}\}$, such that $\{v_{2i-1}, v_{2i}\}$ are roots of the $2$-cmbt of order $i$ $($$2$-$cmbt_i$$)$, for all $1\leq i \leq k-1$.
	\end{enumerate}
\end{definition}

%
%
%\begin{pre}
%	The difference $Tr_2(G)- \frac{Tr(G)}{2}$ or $Tr_2(G)- \frac{Tr(G)+1}{2} $ can be arbitrarily large.
%\end{pre}
%
%\begin{proof}
%	By the Proposition \ref{relation_between_2-transitive_and_transitive} and Theorem \ref{2_transitivity_of_2-cmbt_k}, we have $Tr_2(G)-\frac{Tr(G)}{2}=k-\frac{k}{2}=\frac{k}{2}$ or $Tr_2(G)-\frac{Tr(G)+1}{2}=k-\frac{k+1}{2}=\frac{k-1}{2}$. The difference can be arbitrarily large since $k$ can be made as large as possible. 
%\end{proof}

%and graph having equal transitivity and $2$-transitivity, we define a graph called  and denoted by $2$-$cmbt_k$. The definition is as follows:

\begin{figure}[htbp!]
	\centering
	\includegraphics[scale=0.60]{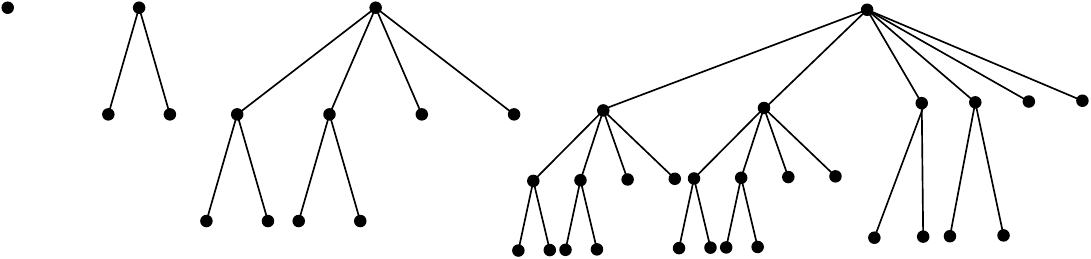}
	\caption{First four $2$-cmbts, $2$-$cmbt_1$, $2$-$cmbt_2$, $2$-$cmbt_3$, and $2$-$cmbt_4$ from the left to right, respectively.}
	\label{fig:first_four_2-cmbt}
\end{figure}

Fig.~\ref{fig:first_four_2-cmbt}, illustrates the $2$-cmbt of orders $1,2,3$ and $4$. In the next lemma, we show that transitivity is the same as $2$-transitivity for such trees.

\begin{lemma}\label{2_transitivity_of_2-cmbt_k}
	Let $T$ be a $2$-$cmbt_k$ ($2$-cmbt of order $k$). Then $Tr_2(T)=Tr(T)=k$.
\end{lemma}
\begin{proof}
	First, we show that $Tr_2(T)=k$. For this, we use induction on the order of $2$-cmbt. For the base case, that is, when $k=1$, by the definition of $2$-cmbt, $T$ is a single vertex graph. Therefore, $Tr_2(T)=1$. Let $T$ be a $2$-$cmbt_k$. By definition, we have $T$ is a rooted tree, whose root has degree $2(k-1)$ and is adjacent to vertices $\{v_1, v_2, \ldots, v_{2(k-1)}\}$, such that $\{v_{2i-1}, v_{2i}\}$ are roots of the $2$-cmbt of order $i$, for all $1\leq i \leq k-1$. Let $c$ be the root of $T$. By the induction hypothesis, $Tr_2(T_{v_{2i-1}})=Tr_2(T_{v_{2i}})=i$, for all $1\leq i \leq k-1$, where $T_v$ is a subtree of $T-\{c\}$, rooted at $v$. Therefore, $Tr_2(T-\{c\})$ is $k-1$. Moreover, we can construct a $2$-transitive partition, say $\pi=\{V_1, V_2, \ldots, V_{k-1}\}$, of $T-\{c\}$ recursively such that  $\{v_{2i-1}, v_{2i}\}\in V_i$ for all $1\leq i \leq k-1$. Now clearly, $\pi'= \pi \cup \{c\}$ forms a $2$-transitive partition of $T$. Therefore, $Tr_2(T)\geq k$. Moreover, by the Proposition \ref{upper_bound_2-transitivity}, we know that $Tr_2(T)\leq \floor{\frac{\Delta(T)}{2}}+1=\floor{\frac{2(k-1)}{2}}+1=k$. Hence, $Tr_2(T)=k$. The algorithm described in \cite{hedetniemi1982linear} shows that $Tr(T)=k$. Therefore, $Tr_2(T)=Tr(T)=k$ for a $2$-$cmbt_k$.
\end{proof}

Next, we characterize the graphs with a small $2$-transitivity.

\begin{proposition}\label{2_transitivity_geq_2}
	Let $G$ be a graph. Then $Tr_2(G)\geq 2$ if and only if $G$ contains $P_3$ as a subgraph.
\end{proposition}
\begin{proof}
	Let $G$ contain $P_3$ as a subgraph. Also, let $v$ be the degree $2$ vertex of $P_3$. Then the partition $\pi =\{V\setminus\{v\}, \{v\}\}$ is a $2$-transitive partition of $G$. Hence, $Tr_2(G)\geq 2$. If $Tr_2(G)\geq 2$, then clearly $G$ contains a $P_3$ as a subgraph, since in any $2$-transitive partition $\pi=\{V_1, V_2, \ldots, V_k\}$, $V_1$ $2$-dominates $V_2$.
\end{proof}

The above proposition also characterizes the connected graphs with $2$-transitivity as $1$. If $Tr_2(G)=1$, then by Proposition \ref{2_transitivity_geq_2}, $G$ does not contain any $P_3$ as a subgraph. Hence, we have the following Corollary.

\begin{corollary}\label{coro_2_transitivity_geq_2}
	For a connected graph $G$, $Tr_2(G)=1$ if and only if $G$ is either $K_1$ or $K_2$.
\end{corollary}
%\begin{proof}
%	Clearly, if $G$ is either $K_1$ or $K_2$, then $Tr_2(G)=1$. Conversely, if $Tr_2(G)=1$, then by Proposition \ref{2_transitivity_geq_2}, $G$ does not contain any $P_3$ as a subgraph. As $G$ is connected, graph $G$ is either $K_1$ or $K_2$.
%\end{proof}

In the following result, we can calculate the $2$-transitivity of a tree $T$ by checking the existence of $2$-cmbt as a subgraph in $T$.

\begin{lemma}
	A tree $T$ has a $2$-transitive partition of size $k$ if and only if $T$ contains a $2$-cmbt of order $k$ as a subgraph. Therefore, $Tr_2(T)=k$ if and only if $T$ contains $2$-cmbt of order $k$ as subgraph for maximum $k$.
\end{lemma}
\begin{proof}
	Since $T$ contains $2$-cmbt of order $k$ as subgraph, Lemma \ref{2_transitivity_of_2-cmbt_k} implies that $Tr_2(G)\geq k$. For the converse, let $T$ have a $2$-transitive partition of size $k$. We show that $T$ contains $2$-cmbt of order $k$ as a subgraph by induction on $k$. The base case for $k=1$ follows from Corollary \ref{coro_2_transitivity_geq_2}, as $T$ is either a single vertex or an edge. In both cases, $T$ contains a $2$-cmbt of order $1$. Let $\pi=\{V_1,V_2, \ldots ,V_k\}$ be a $2$-transitive partition of $T$ of size $k$. Let $F=T\setminus V_1$. Clearly, $Tr_2(F)\geq Tr_2(T)-1=k-1$. Note that $F$ is a forest, and we know that a tree exists, say $T'$, in $F$ such that $Tr_2(T')=Tr_2(F)\geq k-1$. Therefore, $T'$ has a $2$-transitive partition of size $k-1$. By the induction hypothesis, $T'$ contains $2$-cmbt, say $T''$, of order $k-1$ as subgraph. Since $V_1$ $2$-dominates $V_j$ for all $j>1$, for each vertex in $T''$, there exist at least two distinct vertices in $V_1$. So, considering only two neighbours of each vertex in $V_1$ will form a $2$-cmbt tree, a subgraph of $T'$ and $T$ of order $k$. Hence, if $Tr_2(T)=k$, then $T$ contains $2$-cmbt of order $k$ as a subgraph. Let $S$ be the set of vertices of $V_1$ such that each vertex of $T''$ is exactly $2$-dominated by $S$. Note that $|S|=2|V(T'')|$ as $T$ is a tree. From the definition of $2$-cmbt, it follows that subtree induced by $S\cup V(T'')$ is a $2$-cmbt of order $k$. 
	
	Considering the maximum value $k$, we have $Tr_2(T)=k$ if and only if $T$ contains $2$-cmbt of order $k$ as a  subgraph for the maximum $k$.
\end{proof}

Now we show that if the $2$-transitivity of a graph $G$ is constant and $\pi=\{V_1, V_2, \ldots, V_{k-1}, V_k\}$ is a $2$-transitive partition of $G$, then the number of vertices of $V_i$, for $2\leq i\leq k$, can be bound by a constant.

\begin{lemma}\label{size of 2_transitive_partition}
	Let $G$ be a connected graph and $Tr_2(G)=k$, $k\geq 3$. Then there exists a $2$-transitive partition of $G$ of size $k$, say $\pi =\{V_1, V_2, \ldots, V_k\}$ of $G$, such that $|V_k|=1$, $|V_{k-1}|=2$ and $|V_{k-i}|\leq 2\cdot3^{i-1}$, for all $i$, $2\leq i\leq k-2$. 
\end{lemma}
\begin{proof}
	Let $\pi =\{V_1, V_2, \ldots, V_k\}$ be a $2$-transitive partition of $G$ of size $k$, where $k\geq 3$. Now if $|V_k|\geq 2$, then we move all the vertex except one, say vertex $z$, from the set $V_k$ to the first set $V_1$. After this transformation, we have $\pi' =\{V_1'=V_1\cup (V_k-\{z\}), V_2, \ldots, V_{k-1}, V_k'=\{z\}\}$ which is a $2$-transitive partition of $G$ of size $k$ and $|V_k|=1$. Therefore, $\pi'$ is a $2$-transitive of $G$ of size $k$, such that $|V_k|=1$. Let $y_1$ and $y_2$ be two vertices in $V_{k-1}$ that are adjacent to $z$. The existence of these two vertices in $V_{k-1}$ is guaranteed by the fact that $\pi'$ is a $2$-transitive partition of $G$, that is, $V_{k-1}$ $2$-dominates $\{z\}$. Now we move every vertex, except $y_1$ and $y_2$, from $V_{k-1}$ to the set $V'_1$. After doing this transformation, $\pi'$ becomes $\pi''=\{V_1''=V_1'\cup (V_{k-1}-\{y_1, y_2\}), V_2, \ldots, V_{k-2}, V_{k-1}'=\{y_1, y_2\}, V_k'=\{z\}\}$, which is a $2$-transitive partition of $G$ of size $k$ and $|V_{k-1}'|=2$, $|V_k'|=1$. We prove the remaining part of this lemma using induction. Let $\pi=\{V_1, V_2, \ldots, V_k\}$ be a $2$-transitive partition of $G$, such that $|V_k|=1$, $|V_{k-1}|=2$. Consider the base case when $i=2$.  Note that to $2$-dominate $V_k$ and $V_{k-1}$ we need at most $2\cdot(|V_{k-1}|+|V_k|)$ vertices in $V_{k-2}$ and the remaining vertices can be moved to $V_1$ as in the previous argument. Therefore, $|V_{k-2}|\leq 2\cdot(|V_{k-1}|+|V_k|)=2\cdot(2+1)=2\cdot3^{2-1}=2\cdot3^{i-1}$. Assume the induction hypothesis and consider we have a $2$-transitive partition $\pi=\{V_1, V_2, \ldots\, V_k\}$, such that $|V_k|=1$, $|V_{k-1}|=2$ and $|V_{k-j}|\leq 2\cdot3^{j-1}$, $2\leq j\leq i-1$. Now consider the set $V_{k-i}$. Similarly, to $2$-dominate the $i$ sets that are following it, we need at most $2\cdot(|V_{k-i+1}|+ |V_{k-i+2}|+ \ldots+ |V_{k-1}|+|V_k|)$ vertices in $V_{k-i}$. If $V_{k-i}$ contains more than those many vertices, the remaining vertices can be moved to $V_1$, maintaining $2$-transitivity of $G$. Hence, $V_{k-i}\leq 2\cdot(|V_{k-i+1}|+ |V_{k-i+2}|+ \ldots+ |V_{k-1}|+|V_k|)= 2\cdot(2\cdot3^{i-2}+\ldots+2\cdot3+2+1)=2\cdot3^{i-1}$.
\end{proof}
From the above lemma, we have the following corollary:
\begin{corollary}
	The \textsc{M$2$TP} is a fixed parameter tractable when parameterized by the size of the $2$-transitive partition.
\end{corollary}

\begin{remark}
	Note that the converse of Proposition \ref{2-transitivity_equal_transitivity/2} is not true; that is, there are other graphs for which the bound in Proposition \ref{relation_between_2-transitive_and_transitive} is true. For example, consider a split graph $G$, where the complete part has $4$ vertices, the independent part has $3$ vertices, and the degree of each vertex in the independent part is 1. For this graph $G$, $Tr_2(G)=2=\frac{Tr(G)}{2}$ but $4=Tr(G)\not=\Delta(G)+1=4+1$.
\end{remark}

\begin{remark}
	
	Note that by the Propositions \ref{upper_bound_2-transitivity} and \ref{relation_between_2-transitive_and_transitive}, , we have $Tr_2(K_{m, n})= min\{\floor{\frac{m}{2}}+1, \floor{\frac{n}{2}}+1\}$. 
\end{remark}

\begin{remark}
	Note that for transitivity, $Tr(G)=n$ if and only if $G$ is $K_n$ \cite{hedetniemi2018transitivity}. But for $2$-transitivity, $Tr_2(G)=\floor{\frac{n-1}{2}}+1$ does not imply $G$ is a complete graph. For example, consider $G=K_5-M$, where $M$ is the maximum matching of $K_5$. Clearly, $Tr_2(G=K_5-M)=3=\frac{5-1}{2}+1$, but $G$ is not a complete graph.
\end{remark}

\section{NP-completeness}
This section deals with two NP-complete results for M$2$TDP in chordal and bipartite graphs.
\subsection{Chordal graphs}
In this subsection, we show that the \textsc{Maximum $2$-Transitivity Decision Problem} is NP-complete for chordal graphs. This problem is in NP. We prove the NP-completeness of this problem by showing a polynomial-time reduction from \textsc{Proper $3$-Coloring Decision Problem} in graphs having an even number of edges. A \emph{proper $3$-colring} of a graph $G=(V,E)$ is a function $g: V \rightarrow \{1,2,3\}$, such that for every edge $uv \in E$, $g(u)\not= g(v) $. The \textsc{Proper $3$-Coloring Decision Problem} is defined as follows:

\noindent\textsc{\underline{proper $3$-Coloring Decision Problem (P$3$CDP)}}

\noindent\emph{Instance:} A graph $G=(V,E)$

\noindent\emph{Question:} Does there exist a proper $3$-coloring of $V$?

\begin{figure}[htbp!]
	\centering
	\includegraphics[scale=.55]{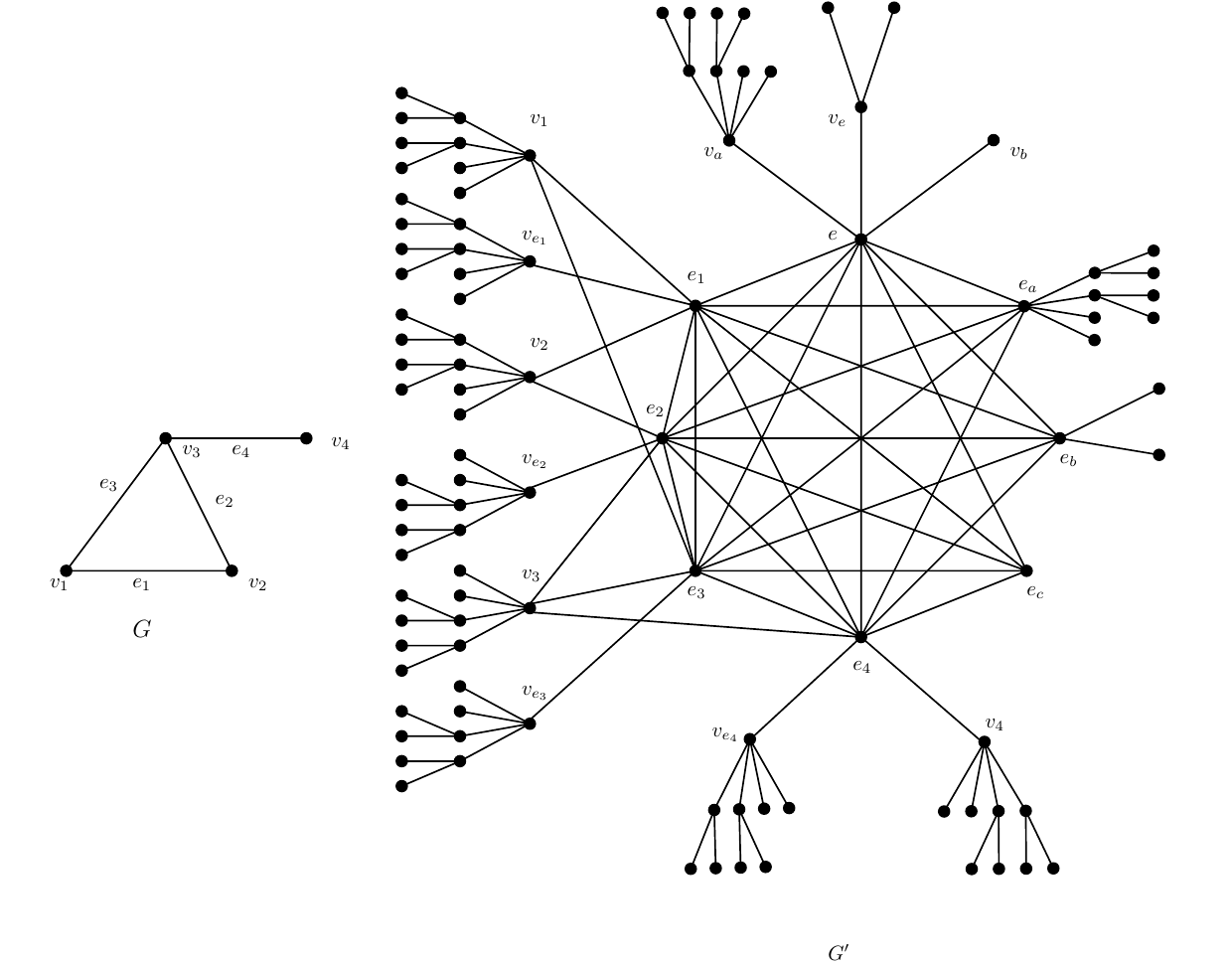}
	\caption{Construction of $G'$}
	\label{fig:chordalnp}
	
\end{figure}

The P$3$CDP is known to be NP-complete for graphs with an even number of edges \cite{garey1990guide}. Given an instance of P$3$CDP, say $G=(V, E)$ with an even number of edges, we construct an instance of M$2$TDP as follows: let $V=\{v_1, v_2, \ldots, v_n\}$ and $E= \{e_1, e_2, \ldots, e_m\}$. For each vertex $v_i\in V$, we consider a $2$-$cmbt_3$ with $v_i$ as its root. Similarly, for each edge $e_i\in E$, we consider another $2$-$cmbt_3$ with $v_{e_i}$ as its root. Further, corresponding to every edge of $E$, we take a vertex $e_j$ in $G'$. Also, we take another vertex $e$ in $G'$. Let $A=\{e_1,e_2,\ldots ,e_m,e\}$. We construct a complete graph with vertex set $A$. Now, for each edge $e_t=v_iv_j\in E$, we join the edges $v_ie_t$, $v_je_t$, and $v_{e_t}e_t$ in $G'$. Next, we consider two $2$-$cmbt_3$ with roots $v_a, e_a$, two $2$-$cmbt_2$ with roots $v_e, e_b$, and two $2$-$cmbt_1$ with roots $v_b,e_c$. Now, we make $e_a, e_b$ and $e_c$ adjacent to every vertex of $A$. Finally, we make $v_a, v_e$, and $v_b$ adjacent to the vertex $e$, and we set $k=\frac{m}{2}+4$. In this construction, we can verify that the graph $G'=(V', E')$ consists of $10m+9n+27$ vertices and $\frac{m^2+29m}{2}+8n+26$ edges. The construction is illustrated in Fig.~\ref{fig:chordalnp}.

Next, we show that $G$ has a proper $3$-coloring if and only if $G'$ has a $2$-transitive partition of size $k$. First, we show the forward direction of this statement in the following lemma.
\begin{lemma}
	If $G=(V,E)$ has a proper $3$-coloring, then $G'=(V', E')$ has a $2$-transitive partition of size $k$.
\end{lemma}
\begin{proof}
	Let $g$ be a proper $3$-coloring of $G$. Based on $g$, let us consider a vertex partition, say $\pi=\{V_1,V_2,\ldots ,V_k\}$ of $G'$ as follows: For each $v_i\in V$, we put the corresponding vertex $v_i$ in $V_q$, if $g(v_i)=q$. For each edge $e_t=v_iv_j\in E$, we put the vertex $v_{e_t}\in V_q$, if $g(v_i)\neq q$ and $g(v_j)\neq q$. Further, we put $v_a, e_a \in V_3$, $v_{e}, e_b \in V_2$, and $v_b, e_c \in V_1$. We put the other vertices of $2$-$cmbt_3$ and $2$-$cmbt_2$ according to the Fig.~\ref{fig:cmbt_coloring}, based on the position of the root vertex in the partition. Finally, for $1\leq j\leq \frac{m}{2}$, we put $e_{2j-1}, e_{2j}\in V_{3+j}$ and $e\in V_{\frac{m}{2}+4}$. Note that $\pi$ is a vertex partition, as $m$ is assumed to be even.
	
	%		\begin{enumerate}
		%			\item  For each $v_i\in V$, we put the corresponding vertex $v_i$ in $V_q$, if $g(v_i)=q$. For each edge $e_t=v_iv_j\in E$, we put the vertex $v_{e_t}\in V_q$, if $g(v_i)\neq q$ and $g(v_j)\neq q$. Further, we put $v_a, e_a \in V_3$, $v_{e}, e_b \in V_2$ and $v_b, e_c \in V_1$. We put the other vertices of $2$-$cmbt_3$ and $2$-$cmbt_2$ according to Figure \ref{fig:cmbt_coloring} based on the position of the root vertex in the partition. Finally, for $1\leq j\leq \frac{m}{2}$, we put $e_{2j-1}, e_{2j}\in V_{3+j}$ and $e\in V_{\frac{m}{2}+4}$. 
		%			
		%				
		%			
		%			
		%			\item For each $v_{e_j}$ vertex corresponding an edge $e_j$ with end points $v_x$ and $v_y$ in $G$, assign $v_{e_j} \in V_l$, where $l= \{1, 2, 3\} \setminus \{g(v_x),g(v_y)\}$.  Put the other vertices of the trees $2$-$cmbt_3$ and $2$-$cmbt_2$ in the partition according to their $v$ vertex and $e_a, e_b$ vertices which are shown in the Figure \ref{fig:cmbt_coloring}.
		%			
		%			\item Let $e_{2j-1}, e_{2j}\in V_{3+j}$, $1\leq j\leq \frac{m}{2}$, and $e\in V_{\frac{m}{2}+4}$. 
		%			
		%	 \end{enumerate}

	\begin{figure}[htbp!]
		\centering
		\includegraphics[scale=0.55]{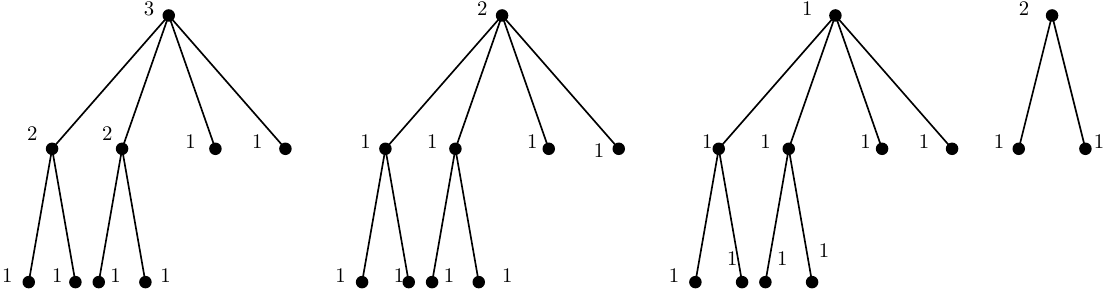}
		\caption{Partition of a $2$-$cmbt_3$ and $2$-$cmbt_2$}
		\label{fig:cmbt_coloring}

	\end{figure}

	Since the vertices of $A$ induce a complete graph, $V_i$ $2$-dominates $V_j$ for $4\leq i<j\leq k$. Also, for each $i=1,2,3$, every vertex of $A$ is adjacent to two vertices of $V_i$. Therefore, for each $i=1,2,3$, $V_i$ $2$-dominates $V_t$, for all $t>3$. Finally, from Fig.~\ref{fig:cmbt_coloring}, it is evident that $V_i$ $2$-dominates $V_j$ for $1\leq i<j\leq 3$. Hence, $\pi$ is a $2$-transitive partition of $G'$ of size $k$. Therefore, if $G$ has a proper $3$-coloring, then $G'$ has a $2$-transitive partition of size $k$.
\end{proof}

%Let $H$ be the complete graph induced by $A$. Since $H$ is a complete graph, then $V_i$ $2$-dominates $V_j$($V_i \doublerightarrow V_j$) for $4\leq i<j\leq k$. Also, each vertex from $A$ is connected with two vertices from $V_i$, for $i=1,2,3$. So, $V_i$, $i=1,2,3$, $2$-dominates $V_t$, for all $t>3$. From the Figure \ref{fig:cmbt_coloring}, it is clear that $V_i$ $2$-dominates $V_j$ for $1\leq i<j\leq 3$. So $\pi$ is a $2$-transitive k-partition of $G'$. Therefore, if $G$ has a proper 3-coloring then $G'$ has a $2$-transitive k-partition.

Next, we show the converse of the statement. For this, we first prove the following claim:

\begin{claim}\label{claim:123 and rest}
	Let  $\pi=\{V_1,V_2,\ldots ,V_k\}$ be a $2$-transitive partition of $G'$ of size $k$ such that $|V_k|=1$. Then the sets $V_4, V_5,\ldots, V_k$ contain only vertices from $A$, and the sets $V_1, V_2,$ and $V_3$ contain only vertices from $V'\setminus A$.
\end{claim}
\begin{proof}
	We divide the proof into two cases:
	
	\begin{case}\label{Chordal_np_two_transitivity_case_1}
		$e\in V_k$
	\end{case}
	First, note that $\{v_a,v_e, v_b\}$ cannot be in $V_p$ for any $p\geq 4$ as their degrees are less than $6$. Therefore, the vertices $\{v_a,v_e,v_b\}$ belong to $V_p$ for some $1\leq p \leq 3$. Since $e\in V_k$, to $2$-dominate $e$, each set in $\{V_1, V_2, \ldots, V_{k-1}\}$ must contain at least two vertices from $N_{G'}(e)= \{e_1, e_2, \ldots, e_m, e_a, e_b, e_c, v_a, v_e, v_b\}$. Since $e$ is adjacent with exactly $2k-2=m+6$ vertices, each $V_i, 1\leq i \leq k-1$ contains exactly two vertices from $\{e_1, e_2, \ldots, e_m, e_a, e_b, e_c, v_a, v_e, v_b\}$. Since $\{v_a,v_e,v_b\}$ belong to $V_p$ for some $1\leq p \leq 3$, it follows that exactly three vertices from $\{e_1, e_2, \ldots, e_m, e_a, e_b, e_c\}$ belong to $V_1, V_2$, and $V_3$. 
	
	Further, note that $e_c$ cannot be in $V_p$ for any $p\geq 4$. Because to $2$-dominate $e_c$, at least six vertices are required in $V_1\cup V_2\cup V_3$ from $\{e_1, e_2, \ldots, e_m, e\}$ which contradicts the fact that exactly one vertex from $\{e_1, e_2, \ldots, e_m, e_a, e_b, e_c\}$ belongs to each of the sets $V_1, V_2$, and $V_3$. For similar reasons, $e_b$ also cannot be in $V_p$ for any $p\geq 4$. Therefore, the vertices $e_c$ and $e_b$ belong to $V_p$ for some $1\leq p \leq 3$. Next, we claim that $e_a$ also cannot be in $V_p$ for any $p\geq 4$. If it happens, then at least two vertices are required in $V_1\cup V_2\cup V_3$ from $\{e_1, e_2, \ldots, e_m\}$ (note that $e\in V_k$). This contradicts the facts that $e_c$ and $e_b$ belong to $V_p$ for some $1\leq p \leq 3$ and exactly one vertex from $\{e_1, e_2, \ldots, e_m, e_a, e_b, e_c\}$ belongs to each of the sets $V_1, V_2$, and $V_3$. Therefore, $e_a$ belongs to $V_p$ for some $1\leq p \leq 3$. Hence, the vertices of $A$ belong to $\{V_4, V_5,\ldots V_k\}$. Note that none of the vertices from $\{v_1, v_2, \ldots, v_n, v_{e_1}, v_{e_2}, \ldots, v_{e_m}\}$ belong to $V_p$ for some $p\geq 4$. Otherwise, two vertices of $A$ must be in $V_3$, as no vertex, except the root, of $2$-$cmbt_3$ can be in $V_3$. But this contradicts that the vertices of $A$ belong to $\{V_4, V_5,\ldots V_k\}$. Since the degree of every other vertices is at most $3$, they cannot belong to $V_p$, $p\geq 4$. Therefore, the sets $V_4, V_5,\ldots V_k$ contain only vertices from $A$ and the sets $V_1, V_2$, and $V_3$ contain only vertices from $V'\setminus A$.

	\begin{case}
		$e\notin V_k$
	\end{case}
	Since the degree of every vertex other than the vertices of $A$ is less than $m+6$ and $e\notin V_k$, without loss of generality, we assume that $V_k=\{e_1\}$, where $e_1$ is the vertex of $G'$ corresponding to the edge $e_1=v_1v_2\in E$. We show that $v_1$ and $v_2$ belong to the first three sets in $\pi$. Let $v_1\in V_l$ and $v_2\in V_t$, where $t\leq l$. If possible, let $l\geq 4$. Since $e_1\in V_k$, to $2$-dominate $e_1$, each set in $\{V_1, V_2, \ldots, V_{k-1}\}$ must contain at least two vertices from $N_{G'}(e_1)=\{e_2, e_3, \ldots, e_m, e, v_1, v_2, v_{e_1}\}$. Since $e_1$ is adjacent to exactly $2k-2=m+6$ vertices, each $V_i$, $1\leq i \leq k-1$, contains exactly two vertices from $N_{G'}(e_1)$. So, if $l\geq 4$, then to $2$-dominates $v_1$, each set in $\{V_3, V_4, \ldots, V_{l-1}\}$ contains exactly two vertices from $\{e_2, e_3, \ldots, e_m\}$. Now, to $2$-dominate $e_1$, each set in $\{V_{l+1}, V_{l+2}, \ldots, V_{k-1}\}$ contains exactly two vertices from $\{e_2, e_3, \ldots, e_m, e, e_a, e_b, e_c\}$. Also, observe that apart from $v_1$, $V_l$ contains exactly one vertex from $\{v_2, e_2, e_3, \ldots, e_m, e, e_a, e_b, e_c\}$. The vertices $\{e, e_a, e_b, e_c\}$ cannot belong to $V_q$ for $q\geq l+1$. Because in that case, $V_l$ would have at least two vertices from $\{e_2, e_3, \ldots, e_m, e, e_a, e_b, e_c\}$, which is a contradiction. Also, the vertices $\{e, e_a, e_b, e_c\}$ cannot belong to $V_p$, $3\leq p \leq l-1$, because we have already seen that $v_1 \in V_l$ and each set in $\{V_3, V_4, \ldots, V_{l-1}\}$ contains exactly two vertices from $\{e_2, e_3, \ldots, e_m\}$ to $2$-dominate $v_1$. Therefore, the vertices $\{e,e_a,e_b,e_c\}$ can be in either $V_1$ or $V_2$, or $V_l$. Hence, each of $\{V_3, V_4, \ldots, V_{l-1},V_{l+1}, \ldots, V_{k-1}\}$ contains exactly two vertices from $\{e_2, \ldots, e_m\}$ only. But the number of vertices in $\{e_2, e_3, \ldots, e_m\}$ is $m-1$, whereas we need $2(k-4)=m$ vertices. Therefore, $l$ cannot be more than $3$. Note that $v_{e_1}$ cannot be in $V_j$ for $j\geq 4$ as its degree is less than $6$. Therefore, the vertices $\{v_1, v_{e_1}, v_2\}$ belong to $V_p$ for $1\leq p \leq 3$. With similar arguments as in Case \ref{Chordal_np_two_transitivity_case_1}, we can say that every vertex of $\{e_a, e_b, e_c\}$ belong $V_p$ for $ 1\leq p \leq 3$, and the vertices of $A$ belong to $\{V_4, V_5,\ldots V_k\}$. We can further claim that the vertices of $\{v_1, v_2, \ldots, v_n, v_{e_1}, v_{e_2}, \ldots, v_{e_m}\}$ belong to $V_p$ for $1\leq p \leq 3$ and the other vertices of $G'$ belong to $V_p$ for $1\leq p \leq 3$. Therefore, the sets $V_4, V_5,\ldots, V_k$ contain only vertices from $A$, and the sets $V_1, V_2$, and $V_3$ contain only vertices from $V'\setminus A$.
\end{proof}

%	\begin{cl}\label{clbp_2}
	%	
	%	The vertex $e_a \in V_3, e_b\in V_2$ and $e_c\in V_1$.
	%	
	%		
	%		\end{cl}
%	
%	
%	\begin{proof}
	%		
	%		Since $e\in V_p, p\geq 4$, to dominate $e$ each of $\{V_1, V_2, V_3\}$ contains exactly two vertices from $\{e_a, e_b, e_c, v_a, v_e, v_b\}$. Clearly, The vertex $e_a \in V_3, e_b\in V_2$ and $e_c\in V_1$.
	%		
	%	\end{proof}

Now we prove $G$ uses $3$-colors, which is a proper coloring in the following lemma.

\begin{lemma}
	If $G'$ has a $2$-transitive partition of size $k$, then $G$ has a proper $3$-coloring.
\end{lemma}
\begin{proof}
	Let $\pi=\{V_1,V_2,\ldots ,V_k\}$ be a $2$-transitive partition of $G'$ of size $k$. By Lemma \ref{size of 2_transitive_partition}, we can assume that $|V_k|=1$. Let us define a coloring of $G$, say $g$, by labelling $v_i$ with color $p$ if its corresponding vertex $v_i$ is in $V_p$. The previous claim ensures that $g$ is a $3$-coloring. Now we show that $g$ is a proper coloring. First note that in $G'$, since $e\in V_p$, with $p\geq 4$, $e_a, e_b$, and $e_c$ must belong to $V_3, V_2$, and $V_1$, respectively. Let $e_t=v_iv_j\in E$, and let its corresponding vertex $e_t$ in $G'$ belong to some set $V_p$ with $p\geq 4$. This implies that the vertices $\{v_i, v_j, v_{e_t}\}$ must belong to different sets from $V_1, V_2$, and $V_3$. Therefore, $g(v_i)\neq g(v_j)$, and hence $g$ is a proper coloring of $G$.
\end{proof}

\subsection{Bipartite graphs}

\begin{figure}[!h]
	\centering
	\includegraphics[scale=0.55]{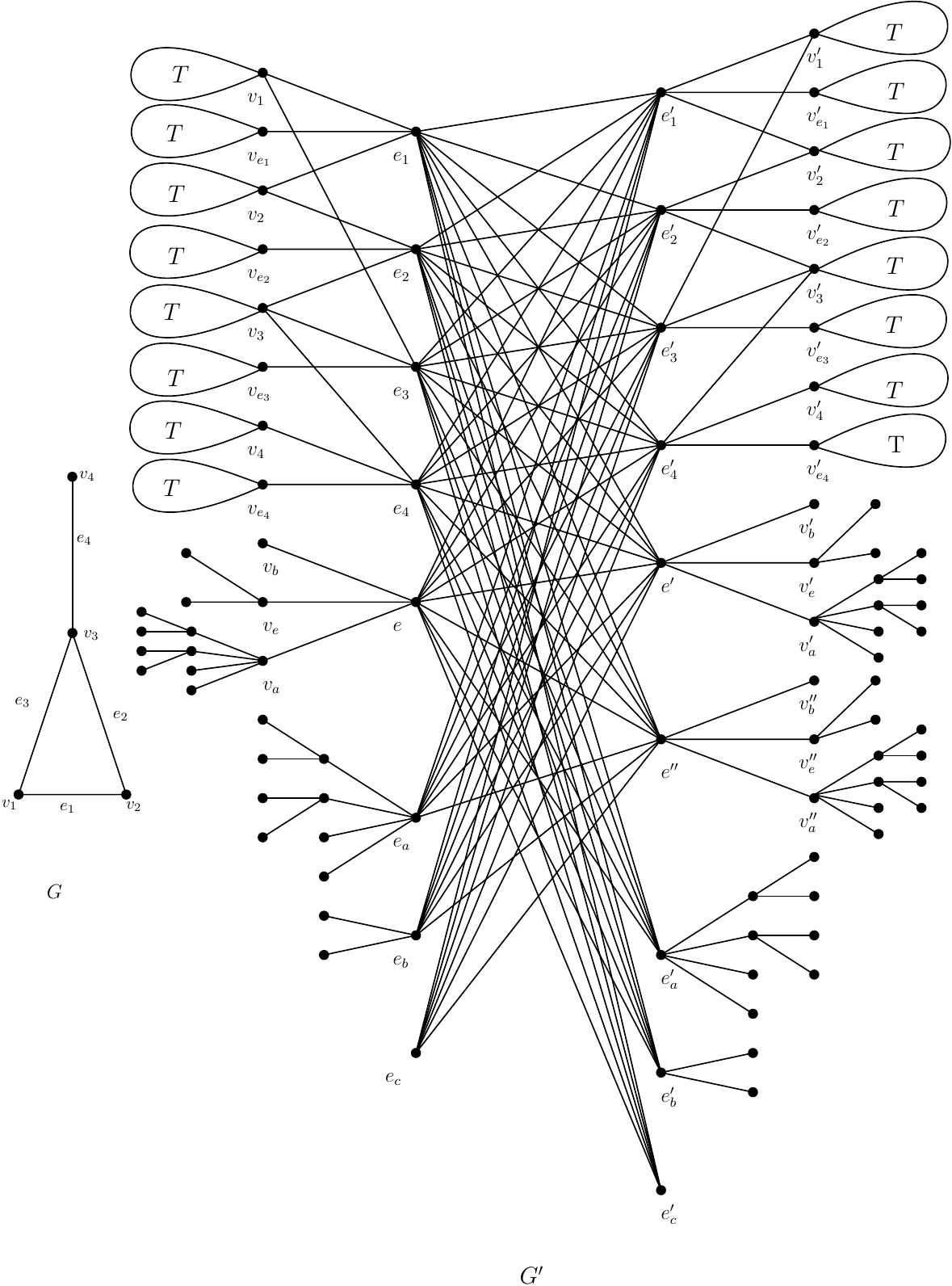}
	\caption{Construction of $G'$, where $T$ is a $2$-$cmbt_3$}
	\label{fig:2_bipartitenp}
	
\end{figure}

In this subsection, we show that the \textsc{Maximum $2$-Transitivity Decision Problem} is NP-complete for a bipartite graph. This problem is in NP. We prove the NP-completeness of this problem by showing a polynomial-time reduction from  \textsc{Proper $3$-Coloring Decision Problem} in graphs having an even number of edges to the M$2$TDP.
% A \emph{proper $3$-colring} of a graph $G=(V,E)$ is a function $g: V \rightarrow \{1,2,3\}$, such that for every edge $uv \in E$, $g(u)\not= g(v) $. The proper \textsc{$3$-Coloring Decision Problem} is defined as follows:
%	
%	\noindent\textsc{\underline{proper $3$-Coloring Decision Problem (P$3$CDP)}}
%	
%	\noindent\emph{Instance:} A graph $G=(V,E)$
%	
%	\noindent\emph{Question:} Does there exist a proper $3$-coloring of $V$?
%	The P$3$CDP is known to be NP-complete for graphs with an even number of edges \cite{garey1990guide}. 
Given an instance of P$3$CDP, say $G=(V, E)$ with an even number of edges, we construct an instance of M$2$TDP as follows: let $V=\{v_1, v_2, \ldots, v_n\}$ and $E= \{e_1, e_2, \ldots, e_m\}$. For each vertex $v_i\in V$, we consider two $2$-$cmbt_3$ with $v_i$ and $v_i'$ as their roots, respectively. Similarly, for each edge $e_i\in E$, we consider another two $2$-$cmbt_3$ with $v_{e_i}$ and $v_{e_i}'$ as their roots, respectively. Further, corresponding to every edge $e_j\in E$, we take two vertices $e_j, e_j'$ in $G'$. Also, we take another three vertices $e, e'$, and $e''$ in $G'$. Let $A=\{e_1, e_2, \ldots , e_m, e\}$ and $B=\{e_1', e_2', \ldots, e_m', e', e'' \}$. We construct a complete bipartite graph with vertex set $A\cup B$. Now, for each edge $e_t=v_iv_j\in E$, we join the edges $v_ie_t$, $v_je_t$, $v_{e_t}e_t$, $v_i'e_t'$, $v_je_t$, and $v_{e_t}'e_t'$ in $G'$. Next, we consider five $2$-$cmbt_3$ with roots $v_a, v_a', v_a'', e_a, e_a'$, five $2$-$cmbt_2$ with roots $v_e, v_e', v_e'', e_b, e_b'$, and five $2$-$cmbt_1$ with roots $v_b, v_b', v_b'', e_c, e_c'$. Now, we make $e_a, e_b$, and $e_c$ adjacent to every vertex of $B$. Also, we make $e_a', e_b'$, and $e_c'$ adjacent to every vertex of $A$. Finally, we make $\{v_a, v_e, v_b\}$, $\{v_a', v_e', v_b'\}$, and $\{v_a'', v_e'', v_b''\}$ adjacent to the vertex $e, e$, and $e''$, respectively, and we set $k=\frac{m}{2}+5$. In this construction, we can verify that the graph $G'=(V', E')$ consists of $18m+18n+68$ vertices and $m^2+31m+16n+70$ edges. The construction is illustrated in Fig.~\ref{fig:2_bipartitenp}.

Next, we show that $G$ has a proper $3$-coloring if and only if $G'$ has a $2$-transitive partition of size $k$. First, we show the forward direction of this statement in the following lemma.

\begin{lemma}
	If $G=(V,E)$ has a proper $3$-coloring, then $G'=(V', E')$ has a $2$-transitive partition of size $k$.
\end{lemma}
\begin{proof}
	Let $g$ be a proper $3$-coloring of $G$. Based on $g$, let us consider a vertex partition, say $\pi=\{V_1,V_2,\ldots ,V_k\}$, of $G'$ as follows: for each $v_i\in V$, we put the corresponding vertex $v_i, v_i'$ in $V_q$, if $g(v_i)=q$. For each edge $e_t=v_iv_j\in E$, we put the vertex $v_{e_t}, v_{e_t}' \in V_q$, if $g(v_i)\neq q$ and $g(v_j)\neq q$. Further, we put $v_a, v_a', v_a'', e_a, e_a' \in V_3$, $v_e, v_e', v_e'', e_b, e_b' \in V_2$, and $v_b, v_b', v_b'', e_c, e_c' \in V_1$. We put the other vertices of $2$-$cmbt_3$ and $2$-$cmbt_2$ according to Fig.~\ref{fig:cmbt_coloring}, based on the position of the root vertex in the partition. Finally, for $1\leq j\leq \frac{m}{2}$, we put $e_{2j-1}, e_{2j}, e_{2j-1}', e_{2j}' \in V_{3+j}$, $e', e''\in V_{\frac{m}{2}+4}$ and $e \in V_{\frac{m}{2}+5}$. Note that $\pi$ is a vertex partition, as $m$ is assumed to be even. Since the vertices of $A\cup B$ induces a complete bipartite graph, $V_i$ $2$-dominates $V_j$ for $4\leq i<j\leq k$. Also, for each $i=1,2,3$, every vertex of $A\cup B$ is adjacent to two vertices of $V_i$. Therefore, for each $i=1,2,3$, $V_i$ $2$-dominates $V_t$, for all $t>3$. Finally, from Fig.~\ref{fig:cmbt_coloring}, it is evident that $V_i$ $2$-dominates $V_j$ for $1\leq i<j\leq 3$. Hence, $\pi$ is a $2$-transitive partition of $G'$ of size $k$. Therefore, if $G$ has a proper $3$-coloring, then $G'$ has a $2$-transitive partition of size $k$.
\end{proof}

Next, we show the converse of the statement. For this, we first prove the following claim:
\begin{claim}\label{claim:2_transitivity_bipartite_np}
	Let  $\pi=\{V_1,V_2,\ldots ,V_k\}$ be a $2$-transitive partition of $G'$ of size $k$ such that $|V_k|=1$ and $|V'_{k-1}|=2$. Then the sets $V_4, V_5,\ldots V_k$ contain only vertices from $A\cup B$, and the sets $V_1, V_2$, and $V_3$ contain only vertices from $V'\setminus (A\cup B)$.
\end{claim}
\begin{proof}
	We divide the proof into four cases:
	
	\begin{case}\label{2BP_case_1}
		$e\in V_k$ and either $e'\in V_{k-1}$ or $e''\in V_{k-1}$
	\end{case}
	First, note that $\{v_a, v_e, v_b\}$ cannot be in $V_p$ for any $p\geq 4$ as their degrees are less than $6$. Similarly, the vertices of $\{v'_a,v'_e,v'_b\}$ and $\{v''_a, v''_e, v''_b \}$ cannot be in $V_p$ for any $p\geq 4$. Therefore, the vertices $\{v_a,v_e,v_b\}$, $\{v'_a,v'_e,v'_b\}$ and $\{v''_a, v''_e, v''_b \}$ belong to $V_p$ for some $1\leq p \leq 3$. Since $e\in V_k$, to $2$-dominate $e$, each set in $\{V_1, V_2, \ldots, V_{k-1}\}$ must contain at least two vertices from $N_{G'}(e)=\{e_1', e_2', \ldots, e_m', e', e'', e_a', e_b', e_c', v_a, v_e, v_b\}$. Since $e$ is adjacent to exactly $2k-2=m+8$ vertices, each $V_i, 1\leq i\leq k-1$ contains exactly two vertices from $N_{G'}(e)$. Since $\{v_a, v_e, v_b\}$ belong to $V_p$ for some $1\leq p \leq 3$, it follows that exactly three vertices from $N_{G'}(e)\setminus \{v_a, v_e, v_b\}$ belong to $V_1, V_2$, and $V_3$. Without loss of generality, let us assume that $e'\in V_{k-1}$. To $2$-dominate $e'$, each set in $\{V_1, V_2, \ldots, V_{k-2}\}$ must contain at least two vertices from $N_{G'}(e') \setminus \{e\}=\{e_1, e_2, \ldots, e_m, e_a, e_b, e_c, v_a', v_e', v_b'\}$ as $e\in V_k$. Since $e'$ is adjacent to exactly $2k-4=m+6$ vertices other than $e$, each $V_i, 1\leq i\leq k-2$ contains exactly two vertices from $N_{G'}(e') \setminus \{e\}$. Similarly, as before, since $\{v_a', v_e', v_b'\}$ belong to $V_p$ for some $1\leq p \leq 3$, it follows that exactly three vertices from $N_{G'}(e')\setminus \{e, v'_a, v'_e, v'_b\}$ belong to $V_1, V_2$, and $V_3$. Further, note that $e_c'$ cannot be in $V_p$ for any $p\geq 4$. Because to $2$-dominate $e_c$, at least six vertices are required in $V_1\cup V_2\cup V_3$ from $\{e_1, e_2, \ldots, e_m, e\}$, which contradicts the fact that exactly three vertices from $N_{G'}(e')\setminus \{e, v_a, v_e, v_b\}$ belong to $V_1, V_2$, and $V_3$. For similar reasons, $e_b'$ also cannot be in $V_p$ for any $p\geq 4$. Therefore, the vertices $e_c'$ and $e_b'$ belong to $V_p$ for some $1\leq p \leq 3$. Similarly, we can argue that the vertices $e_c$ and $e_b$ belong to $V_p$ for some $1\leq p \leq 3$. 
	
	Next, we show that $e_a'$ also cannot be in $V_p$ for any $p\geq 4$. If it happens, then at least two vertices are required in $V_1\cup V_2\cup V_3$ from $\{e_1, e_2, \ldots, e_m, e\}$ (note that $e\in V_k$). This contradicts the facts that $e_c$ and $e_b$ belong to $V_p$ for some $1\leq p \leq 3$ and exactly three vertices from $N_{G'}(e')\setminus \{e, v'_a, v'_e, v'_b\}$ belong to $V_1, V_2$, and $V_3$. Therefore, $e_a'$ belongs to $V_p$ for some $1\leq p \leq 3$. Similarly, we can show that $e_a$ belongs to $V_p, 1\leq p \leq 3$. Hence, the vertices of $A\cup B$ belong to $\{V_4, V_5,\ldots V_k\}$. Note that none of the vertices from $\{v_1, v_2, \ldots, v_n, v_{e_1}, v_{e_2}, \ldots, v_{e_m}\}$ belong to $V_p$ for some $p\geq 4$. Otherwise, two vertices of $B$ must be in $V_3$, as no vertex, except the root, of $2$-$cmbt_3$ can be in $V_3$. But this contradicts the fact that the vertices of $B$ belong to $\{V_4, V_5,\ldots V_k\}$. Similarly, no vertices from $\{v_1', v_2', \ldots, v_n', v_{e_1}', v_{e_2}', \ldots, v_{e_m}'\}$ belong to $V_p$ for some $p\geq 4$. Since the degree of every other vertices is at most $3$, they cannot belong to $V_p$, $p\geq 4$. Therefore, the sets $V_4, V_5,\ldots V_k$ contain only vertices from $A\cup B$, and the sets $V_1, V_2$, and $V_3$ contain only vertices from $V'\setminus (A\cup B)$.
	
	\begin{case}\label{2BP_case_2}
		$e\in V_k$ and $e', e''\notin V_{k-1}$
	\end{case}
	
	First, note that $\{v_a, v_e, v_b\}$ cannot be in $V_p$ for any $p\geq 4$ as their degrees are less than $6$. Since $e\in V_k$, to $2$-dominate $e$, each set in $\{V_1, V_2, \ldots, V_{k-1}\}$ must contain at least two vertices from $N_{G'}(e)$, which is equal to the set $\{e_1', e_2', \ldots, e_m', e', e'', e_a', e_b', e_c', v_a, v_e, v_b\}$. Since $e$ is adjacent to exactly $2k-2=m+8$ vertices, each $V_i, 1\leq i\leq k-1$ contains exactly two vertices from $N_{G'}(e)$. Since $\{v_a, v_e, v_b\}$ belong to $V_p$ for some $1\leq p \leq 3$, it follows that exactly three vertices from $N_{G'}(e)\setminus \{v_a, v_e, v_b\}$ belong to $V_1, V_2$, and $V_3$. Since the degree of $\{e_a', e_b', e_c'\}$ is less than $m+6$, these vertices cannot belong to $V_{k-1}$. Again, we know that $e', e''\notin V_{k-1}$, then to $2$-dominate $e$, there exist two vertices in $\{e'_1, e'_2, \ldots, e'_m\} $, say $e'_i, e'_j$, such that $e'_i, e'_j \in V_{k-1}$. Let $e_j$ be incident to $v_r$ and $v_{s}$ in $G$. Now, to $2$-dominate $e'_j$, each set in $\{V_1, V_2, \ldots, V_{k-2}\}$ must contain at least two vertices from $N_{G'}(e'_j) \setminus \{e\}=\{e_1, e_2, \ldots, e_m, e_a, e_b, e_c, v_r', v_{e_j}', v_{s}'\}$ as $e\in V_k$. Since $e'_j$ is adjacent to exactly $2k-4=m+6$ vertices other than $e$, each $V_i, 1\leq i\leq k-2$ contains exactly two vertices from $N_{G'}(e'_j) \setminus \{e\}$. Note that the vertex $v_{e_j}'$ cannot be in $V_p$, for some $p\geq 4$, as the degree of $v_{e_j}'$ is $5$. 
	
	Further, note that $e_c$ cannot be in $V_p$ for any $p\geq 4$. Because to $2$-dominate $e_c$, at least six vertices are required in $V_1\cup V_2\cup V_3$ from $\{e_1', e_2', \ldots, e_m', e', e''\}$ which contradicts the fact that exactly three vertices from $N_{G'}(e)\setminus \{v_a, v_e, v_b\}$ belong to $V_1, V_2$ and $V_3$. For similar reasons, $e_b$ also cannot be in $V_p$ for any $p\geq 4$. Therefore, the vertices $e_c$ and $e_b$ belong to $V_p$ for some $1\leq p \leq 3$. Since $\{e_b, e_c, v_{e_j}'\}$ belong to $V_p$ for some $1\leq p \leq 3$, it follows that exactly three vertices from $N_{G'}(e'_j)\setminus \{e, e_b, e_c, v_{e_j}'\}$ belong to $V_1, V_2$, and $V_3$. Moreover, $e_c'$ cannot be in $V_p$ for any $p\geq 4$. Because to $2$-dominate $e_c'$, at least six vertices are required in $V_1\cup V_2\cup V_3$ from $\{e_1, e_2, \ldots, e_m, e\}$ which contradicts the fact that exactly three vertices from $N_{G'}(e'_j)\setminus \{e, e_b, e_c, v_{e_j}'\}$ belong to $V_1, V_2$ and $V_3$. For similar reasons, $e_b'$ also cannot be in $V_p$ for any $p\geq 4$. Therefore, the vertices $e_c'$ and $e_b'$ belong to $V_p$ for some $1\leq p \leq 3$. 
	
	Next, we show that $e_a$ also cannot be in $V_p$ for any $p\geq 4$. If it happens, then at least two vertices are required in $V_1\cup V_2\cup V_3$ from $\{e_1', e_2', \ldots, e_m', e', e''\}$. This contradicts the facts that $e_c'$ and $e_b'$ belong to $V_p$ for some $1\leq p \leq 3$ and exactly three vertices from $N_{G'}(e)\setminus \{v_a, v_e, v_b\}$ belong to $V_1, V_2$, and $V_3$. Therefore, $e_a$ belongs to $V_p$ for some $1\leq p \leq 3$. Now, $v_r'$ cannot be in $V_p$, for some $p\geq 4$. If it happens, at least two vertices are required in $V_1\cup V_2\cup V_3$ from $\{e_1', e_2', \ldots, e_m'\}$. This contradicts the facts that $e_c'$ and $e_b'$ belong to $V_p$ for some $1\leq p \leq 3$ and exactly three vertices from $N_{G'}(e)\setminus \{v_a, v_e, v_b\}$ belong to $V_1, V_2$, and $V_3$. Therefore, $v_j'$ belongs to $V_p$ for some $1\leq p \leq 3$. Similarly, the vertex $v'_{s}$ belongs to $V_p$ for some $1\leq p \leq 3$. Therefore, the vertices $\{e_a, e_b, e_c, v_r', v_{e_j}', v'_{s}\}$ belong to $V_p$ for some $1\leq p \leq 3$. We show that $e_a'$ also cannot be in $V_p$ for any $p\geq 4$. If it happens, then at least two vertices are required in $V_1\cup V_2\cup V_3$ from $\{e_1, e_2, \ldots, e_m\}$ (note that $e\in V_k$). This contradicts the facts that vertices of $\{e_a, e_b, e_c, v_j', v_{e_j}', v_{j+1}\}$ belong to $V_p$ for some $1\leq p \leq 3$, and each $V_i, 1\leq i\leq k-2$ contains exactly two vertices from $N_{G'}(e'_j) \setminus \{e\}$. Therefore, $e_a'$ belongs to $V_p$ for some $1\leq p \leq 3$. Hence, the vertices of $A\cup B$ belong to $\{V_4, V_5,\ldots V_k\}$. Now, using similar arguments as in Case \ref{2BP_case_1}, we can say that $\{V_4, V_5,\ldots V_k\}$ contain only vertices from $A\cup B$ and $\{V_1, V_2, V_3\}$ contain only vertices from $V'\setminus (A\cup B)$.
	
	\begin{case}\label{2BP_case_3}
		$e\notin V_k$ and either $e'\in V_{k-1}$ or $e''\in V_{k-1}$
	\end{case}
	
	Now, if a vertex $v\in V_k$, then $deg(v)\geq m+8$. Only vertices from $\{e_1, e_2, \ldots, e_m, e\}$ have $deg(v)\geq m+8$. Since $e\notin V_k$, without loss of generality, let us assume that $e_1\in V_k$ and $e_1\in E$ be the edge between $v_1$ and  $v_2$ in $G$. Also, assume $e'\in V_{k-1}$. Note that $\{v_a', v_e', v_b'\}$ cannot be in $V_p$ for any $p\geq 4$ as their degrees are less than $6$. To $2$-dominate $e'$, each set in $\{V_1, V_2, \ldots, V_{k-2}\}$ must contain at least two vertices from $N_{G'}(e') \setminus \{e_1\}=\{e_1, e_2, \ldots, e_m, e_a, e_b, e_c, v_a', v_e', v_b'\}$ as $e_1\in V_k$. Since $e'$ is adjacent to exactly $2k-4=m+6$ vertices other than $e_1$, each $V_i, 1\leq i\leq k-2$ contains exactly two vertices from $N_{G'}(e') \setminus \{e_1\}$. Since $\{v_a', v_e', v_b'\}$ belong to $V_p$ for some $1\leq p \leq 3$, it follows that exactly three vertices from $N_{G'}(e')\setminus \{e_1, v'_a, v'_e, v'_b\}$ belong to $V_1, V_2$, and $V_3$. Since $e_1\in V_k$, to $2$-dominate $e$, each set in $\{V_1, V_2, \ldots, V_{k-1}\}$ must contain at least two vertices from $N_{G'}(e_1)=\{e_1', e_2', \ldots, e_m', e', e'', e_a', e_b', e_c', v_1, v_{e_1}, v_2\}$. Since $e_1$ is adjacent to exactly $2k-2=m+8$ vertices, each $V_i, 1\leq i\leq k-1$ contains exactly two vertices from $N_{G'}(e)$. Note that the vertex $v_{e_1}$ cannot be in $V_p$, for some $p\geq 4$, as the degree of $v_{e_1}$ is $5$. 
	
	Further, note that $e_c'$ cannot be in $V_p$ for any $p\geq 4$. Because to $2$-dominate $e_c'$ at least six vertices are required in $V_1\cup V_2\cup V_3$ from $\{e_1, e_2, \ldots, e_m, e\}$ which contradicts the fact that exactly three vertices from $N_{G'}(e')\setminus \{e_1, v_a', v_e', v_b'\}$ belong to $V_1, V_2$, and $V_3$. For similar reasons, $e_b'$ also cannot be in $V_p$ for any $p\geq 4$. Therefore, the vertices $e_c'$ and $e_b'$ belong to $V_p$ for some $1\leq p \leq 3$. Since $\{e_b', e_c', v_{e_1}\}$ belong to $V_p$ for some $1\leq p \leq 3$, it follows that exactly three vertices from $N_{G'}(e_1)\setminus \{e_b', e_c', v_{e_1}\}$ belong to $V_1, V_2$, and $V_3$. Moreover, $e_c$ cannot be in $V_p$ for any $p\geq 4$. Because to $2$-dominate $e_c$ at least six vertices are required in $V_1\cup V_2\cup V_3$ from $\{e_1', e_2', \ldots, e_m', e', e''\}$ which contradicts the fact that exactly exactly three vertices from $N_{G'}(e_1)\setminus \{e_b', e_c', v_{e_1}\}$ belong to $V_1, V_2$, and $V_3$. For similar reasons, $e_b$ also cannot be in $V_p$ for any $p\geq 4$. Therefore, the vertices $e_c$ and $e_b$ belong to $V_p$ for some $1\leq p \leq 3$. 
	
	Next, we show that $e_a'$ also cannot be in $V_p$ for any $p\geq 4$. If it happens, at least two vertices are required in $V_1\cup V_2\cup V_3$ from $\{e_1, e_2, \ldots, e_m, e\}$. This contradicts the facts that $e_c$ and $e_b$ belong to $V_p$ for some $1\leq p \leq 3$ and exactly three vertices from $N_{G'}(e')\setminus \{e_1, v_a', v_e', v_b'\}$ belong to $V_1, V_2$, and $V_3$. Therefore, $e_a'$ belongs to $V_p$ for some $1\leq p \leq 3$. Now, $v_1$ cannot be in $V_p$, for some $p\geq 4$. If it happens, at least two vertices are required in $V_1\cup V_2\cup V_3$ from $\{e_1, e_2, \ldots, e_m\}$. This contradicts the facts that $e_c$ and $e_b$ belong to $V_p$ for some $1\leq p \leq 3$ and exactly three vertices from $N_{G'}(e')\setminus \{e_1, v_a', v_e', v_b'\}$ belong to $V_1, V_2$, and $V_3$. Therefore, $v_1$ belongs to $V_p$ for some $1\leq p \leq 3$. Similarly, the vertex $v_{2}$ belongs to $V_p$ for some $1\leq p \leq 3$. Therefore, the vertices $\{e_a', e_b', e_c', v_1, v_{e_1}, v_2\}$ belong to $V_p$ for some $1\leq p \leq 3$. We show that $e_a$ also cannot be in $V_p$ for any $p\geq 4$. If it happens, then at least two vertices are required in $V_1\cup V_2\cup V_3$ from $\{e_1', e_2', \ldots, e_m', e', e''\}$ (note that $e'\in V_k$). This contradicts the facts that vertices of $\{e_a', e_b', e_c', v_1, v_{e_1}, v_2\}$ belong to $V_p$ for some $1\leq p \leq 3$, and each $V_i, 1\leq i\leq k-1$ contains exactly two vertices from $N_{G'}(e_1)$. Therefore, $e_a$ belongs to $V_p$ for some $1\leq p \leq 3$. Now by using similar arguments as in Case \ref{2BP_case_1}, we can say that $\{V_4, V_5,\ldots V_k\}$ contain only vertices from $A\cup B$ and $\{V_1, V_2, V_3\}$ contain only vertices from $V'\setminus (A\cup B)$.

	\begin{case}\label{2BP_case_4}
		$e\notin V_k$ and  $e', e''\notin V_{k-1}$
	\end{case}
	
	Note that if a vertex $v$ belongs to $V_k$, then $deg_{G'}(v)\geq m+8$. Only the vertices in $\{e_1, e_2, \ldots, e_m, e\}$ have a degree greater or equal to $m+8$. Since $e\notin V_k$, without loss of generality, let us assume that $e_1\in V_k$ and $e_1$ is the edge between $v_1$ and $v_2$ in $G$. To $2$-dominate $e_1$, each set in $\{V_1, V_2, \ldots, V_{k-1}\}$ must contain at least two vertices from $N_{G'}(e_1)$, that is, equal to $\{e_1', e_2', \ldots, e_m', e', e'', e_a', e_b', e_c', v_1, v_{e_1}, v_2\}$. Since $e_1$ is adjacent with exactly $2k-2=m+8$ vertices, each $V_i$ contains exactly two vertices from $N_{G'}(e_1)$ for $1\leq i\leq k-1$. Since the degree of each vertex in $\{e_a', e_b', e_c', v_1, v_2, v_{e_1}\}$ is less than $m+6$, these vertices cannot be in $V_{k-1}$. Again, we know that $e', e''\notin V_{k-1}$. Therefore, to $2$-dominate $e_1$, there exist two vertices in $\{e'_1, e'_2, \ldots, e'_m\} $, say $e'_i, e'_j$, such that $e'_i, e'_j \in V_{k-1}$. Let $e_j$ be incident to $v_r$ and $v_s$ in $G$. Now, to $2$-dominate $e'_j$, each set in $\{V_1, V_2, \ldots, V_{k-2}\}$ must contain at least two vertices from $N_{G'}(e'_j) \setminus \{e_1\}=\{e_2, \ldots, e_m, e, e_a, e_b, e_c, v_r', v_{e_j}', v_{s}'\}$ as $e_1\in V_k$. Since $e'_j$ is adjacent with exactly $2k-4=m+6$ vertices other than $e_1$, each $V_i, 1\leq i\leq k-2$ contains exactly two vertices from $N_{G'}(e'_j) \setminus \{e_1\}$.
	
	We show that $v_1$ and $v_2$ belong to the first three sets in $\pi$. Let $v_1\in V_l$ and $v_2\in V_t$, where $t\leq l$. If possible, let $l\geq 4$. In that case, to $2$-dominate $v_1$, each set in $\{V_3, V_4, \ldots, V_{l-1}\}$ contains exactly two vertices from $\{e_2, e_3, \ldots, e_m\}$. Therefore, vertices from $\{e, e_a, e_b, e_c, v_r', v_s', v_{e_j}'\}$ cannot be in $\{V_3, V_4, \ldots, V_{l-1}\}$.
	
	Now, we argue that the vertices $\{e, e_a, e_b, e_c, v_r', v_s', v_{e_j}'\}$ cannot be in $V_q$ for $q\geq l+1$. We show this fact for the vertex $e$. First observe that apart from $v_1$, $V_l$ contains exactly one vertex from $N_{G'}(e_1)\setminus\{v_1\}$, that is,  $\{v_{e_1}, v_2, e_1', e_2', e_3', \ldots, e_m', e',e'', e_a', e_b', e_c'\}$. Also, $\{v_a, v_e, v_b\}$ cannot be in $V_p$ for any $p\geq 4$ as their degrees are less than $6$. Now, if $e\in V_p$ for some $p\geq l+1$, then to $2$-dominate $e$, in $V_l$ we require at least two vertices from $N_{G'}(e)\setminus \{v_a, v_e, v_b\}$, which is the same as $\{e_1', e_2', e_3', \ldots, e_m', e', e'', e_a', e_b', e_c'\}$. Hence, we have a contradiction as we have seen that $V_l$ contains exactly one vertex from $N_{G'}(e_1)\setminus\{v_1\}=\{v_{e_1}, v_2, e_1', e_2', e_3', \ldots, e_m', e',e'', e_a', e_b', e_c'\}$. Using similar arguments, we can show that the vertices $\{e_a, e_b, e_c, v_r', v_s', v_{e_j}'\}$ cannot be in $V_p$ for some $p\geq l+1$. Therefore, the vertices $\{e, e_a, e_b, e_c, v_r', v_s', v_{e_j}'\}$ can be in either $V_1$, $V_2$, or $V_l$. It follows that to $2$-dominate $e_j'$, each set in $\{V_{l+1}, V_{l+2}, \ldots, V_{k-2}\}$ contains exactly two vertices from $\{e_2, e_3, \ldots, e_m\}$. Hence, each of $\{V_3, V_4, \ldots, V_{l-1},V_{l+1}, \ldots, V_{k-2}\}$ contains exactly two vertices from $\{e_2, \ldots, e_m\}$ only. But the number of vertices in $\{e_2, e_3, \ldots, e_m\}$ is $m-1$, whereas we need $2(k-5)=m$ vertices. Therefore, $l$ cannot be more than $3$.
	
	Now, $v_{e_1}$ cannot be in $V_p$ for $p\geq 4$ as its degree is $5$. Therefore, the vertices $\{v_1, v_{e_1}, v_2\}$ belong to $V_p$ for $1\leq p \leq 3$. Therefore, we are in a similar situation as in Case \ref{2BP_case_2}, where a vertex from $A$ is in $V_k$, its neighbours outside $A\cup B$ are in $V_1\cup V_2\cup V_3$ and $e', e''\notin V_{k-1}$. By using similar arguments as in Case \ref{2BP_case_2}, we can say that $\{V_4, V_5,\ldots V_k\}$ contain only vertices from $A\cup B$ and $\{V_1, V_2, V_3\}$ contain only vertices from $V'\setminus (A\cup B)$.
\end{proof}
Now we prove $G$ uses $3$-colors, which is a proper coloring in the following lemma.

\begin{lemma}
	If $G'$ has a $2$-transitive partition of size $k$, then $G$ has a proper $3$-coloring.
\end{lemma}

\begin{proof}
	Let $\pi=\{V_1,V_2,\ldots,V_k\}$ be a $2$-transitive partition of $G'$ of size $k$. By Lemma \ref{size of 2_transitive_partition}, we can assume that $|V_k|=1$ and $|V_{k-1}|=2$. Let us define a coloring of $G$, say $g$, by labelling $v_i$ with color $p$ if its corresponding vertex $v_i$ is in $V_p$. The previous claim ensures that $g$ is a $3$-coloring. Now we show that $g$ is a proper coloring. First note that in $G'$, since $e\in V_p$, with $p\geq 4$, $e_a', e_b'$, and $e_c'$ must belong to $V_3, V_2$, and $V_1$, respectively. Let $e_t=v_iv_j\in E$ and its corresponding vertex $e_t$ in $G'$ belong to some set $V_p$ with $p\geq 4$. This implies that the vertices $\{v_i, v_j, v_{e_t}\}$ must belong to different sets from $V_1, V_2$, and $V_3$. Therefore, $g(v_i)\neq g(v_j)$, and hence $g$ is a proper $3$-coloring of $G$.
\end{proof}

\section{Linear-time algorithms}
In this section, we present three linear-time algorithms to find $2$-transitivity for trees, split graphs and bipartite chain graphs.

\subsection{Trees}
In this subsection, we design a linear-time algorithm for finding the $2$-transitivity of a given tree $T=(V, E)$. Our algorithm is similar to the algorithm for finding the Grundy number of an input tree presented in \cite{hedetniemi1982linear}. First, we give a comprehensive description of our proposed algorithm.

\subsubsection{Description of the algorithm}
Let $T^c$ denote a rooted tree rooted at a vertex $c$, and $T_v^c$ denote the subtree of $T^c$ rooted at a vertex $v$. With a small abuse of notation, we use $T^c$ to denote both the rooted tree and the underlying tree. To find the $2$-transitivity of $T=(V,E)$, we first define the \emph{$2$-transitive number} of a vertex $v$ in $T$. The $2$-transitive number of a vertex $v$ in $T$ is the maximum integer $p$ such that $v\in V_p$ in a $2$-transitive partition $\pi=\{V_1, V_2, \ldots, V_k\}$, where the maximum is taken over all $2$-transitive partitions of $T$. We denote the $2$-transitive number of a vertex $v$ in $T$ by $t_2(v, T)$. Note that the $2$-transitivity of $T$ is the maximum $2$-transitive number that a vertex can have; that is, $Tr_2(T)=\max\limits_{v\in V}\{t_2(v, T)\}$. Therefore, our goal is to find the $2$-transitive number of every vertex in the tree. To this end, we define another parameter, namely the \emph{rooted $2$-transitive number}. The \emph{rooted $2$-transitive number} of $v$ in $T^c$ is the $2$-transitive number of $v$ in the tree $T_v^c$, and it is denoted by $t_2^r(v, T^c)$. Therefore,  $t_2^r(v, T^c)= t_2(v, T_v^c)$. Note that the value of the rooted $2$-transitive number of a vertex is dependent on the rooted tree, whereas the $2$-transitive number is independent of the rooted tree. Also, for the root vertex $c$, $t_2^r(c, T^c)= t_2(c,T)$. We recursively compute the rooted $2$-transitive number of the vertices of $T^c$ in a bottom-up approach. First, we consider a vertex ordering $\sigma$, the reverse of the BFS ordering of $T^c$. For a leaf vertex $c_i$, we set $t^r_2(c_i, T^c)=1$. For a non-leaf vertex $c_i$, we call the function \textsc{$2$-Transitive\_Number$()$}, which takes the rooted $2$-transitive number of children of $c_i$ in $T^c$ as input and returns the rooted $2$-transitive number of $c_i$ in $T^c$. At the end of the bottom-up approach, we have the rooted $2$-transitive number of $c_i$ in $T^c$, that is, $t_2^r(c, T^c)$, which is the same as the $2$-transitive number of $c$ in $T$, that is, $t_2(c, T)$. After the bottom-up approach, we have the $2$-transitive number of the root vertex $c$ and the rooted $2$-transitive number of every other vertices in $T^c$.

\begin{algorithm}[!h]
	\footnotesize
	\caption{\textsc{$2$-Transitivity(T)} } \label{Algo:2-trasitivity(T)}
	
	\textbf{Input:} A tree $T=(V, E)$.
	
	\textbf{Output:} $2$-transitivity of $T$.
	
	\begin{algorithmic}[1]

		%\State Construct a rooted tree $T^c$, rooted at a vertex $c$.
		
		\State  Let $\sigma=(c_1, c_2, \ldots, c_k=c)$ be the reverse BFS ordering of the  vertices of $T^c$, rooted at a vertex $c$. 
		
		\ForAll {$c_i$ in $\sigma$}
		
		\If {$c_i$ is a leaf }
		
		\State $t^r_2(c_i, T^c)=1$.

		\Else
		
		\State $t^r_2(c_i, T^c)$ = \textsc{$2$-Transitive\_Number}$(t^r_2(c_{i_1}, T^c), \ldots , t^r_2(c_{i_k}, T^c) )$.
		
		~~~~~~~~~~~~~ /*where $c_{i_1}, \ldots, c_{i_k}$ be the children of $c_i$ and $t^r_2(c_{i_1}, T^c)\leq  \ldots \leq t^r_2(c_{i_k}, T^c)$*/
		
		\EndIf
		
		\EndFor
		
		%\State $z=t_2(c, T)$.
		
		%\State Let $c_{1}, c_{2}, \ldots, c_{p}$ be the children of $c$ in $T^c$ and $t^r_2(c_{1}, T^c)\leq t^r_2(c_{2}, T^c)\leq \ldots \leq t^r_2(c_{p}, T^c)$.
		
		\State \textsc{Mark\_Required}$(t_2(c, T), t^r_2(c_{1}, T^c),  \ldots, t^r_2(c_{k}, T^c) )$
		
		~~~~~~~~~~~~~~/*where $c_{1}, \ldots, c_{k}$ be the children of $c$ and $t^r_2(c_{1}, T^c)\leq  \ldots \leq t^r_2(c_{k}, T^c)$*/
		
		%\State Let $\sigma'$ be the BFS ordering of the vertices of $T^c$
		
		\ForAll{$c_i\in \sigma'\setminus\{c\}$ } ~~~~~~~~~~~~~~~~~~~~~~~~~~~~~~~~~~~~~~~~~~~~~~~~~~~~~/*where $\sigma'$ is the BFS ordering of $T^c$*/
		
		\If {$R(c_i)=0$}
		
		\State $t^r_2(p(c_i), T^{c_i}) = t_2(c_i, T)$
		
		\Else
		
		\State $t^r_2(p(c_i), T^{c_i}) = t_2(c_i, T) -1$
		
		\EndIf

		%\State  Let $c_{i_1}, c_{i_2}, \ldots, c_{i_p}$ be the neighbour of $c_i$ in $T$, where $c_{i_1}, c_{i_2}, \ldots, c_{i_{p-1}}$ are the children of $c_i$ in $T^c$ and $c_{i_p}$ is the parent of $c_i$ in $T^c$.

		%\State Calculate $t^r_2(c_{i_j}, T^v_i)$ for all $j\in \{1, 2, \ldots, q\}$
		
		%\State $t^r_2(c_{i_j}, T^{c_i})=t^r_2(c_{i_j}, T^{c})$ for all $i_1\leq s\leq i_{p-1}$ and $t^r_2(c_{i_p}, T^{c_i})=t_2(c, T)$, if $c_i$ is not ``Required'', otherwise $t^r_2(c_{i_p}, T^{c_i})=t_2(c, T)-1$
		
		\State $t_2(c_i, T) =$\textsc{$2$-Transitive\_Number}$(t^r_2(c_{i_1}, T^{c_i}), \ldots , t^r_2(c_{i_k}, T^{c_i}) )$.

		\State \textsc{Mark\_Required}$(t_2(c_i, T), t^r_2(c_{i_1}, T^{c_i}), \ldots , t^r_2(c_{i_k}, T^{c_i}) )$

		~~~~~~~~~~~~~~~~~~~~~~~~~/*where $c_{i_1}, \ldots, c_{i_k}$ be the neighbours of $c_i$ and $t^r_2(c_{i_1}, T^c)\leq  \ldots \leq t^r_2(c_{i_k}, T^c)$*/
		
		~~~~~~~~~~~~~~~~~~~~~~~~~~~~~~~~~~~~~~~~/* Note that $t^r_2(c_{i_j}, T^{c_i})=t^r_2(c_{i_j}, T^{c})$ if $c_{i_j}$ is a child of $c_i$ in $T^c$*/
		
		\EndFor

		\State $Tr_2(T)=\max\limits_{x\in V}\{t_2(x, T)\}$.

	\end{algorithmic}
	
\end{algorithm}

Next, we compute the $2$-transitive number of every other vertex. For a vertex $c_i$, other than $c$, we compute the $2$-transitive number using \textsc{$2$-Transitive\_Number$()$}, which takes the rooted $2$-transitive number of children of $c_i$ in $T^{c_i}$ as input. Let $y$ be the parent of $c_i$ in $T^c$. Note that, except for $y$, the rooted $2$-transitive number of children of $c_i$ in $T^{c_i}$ is the same as the rooted $2$-transitive number in $T^c$. We only need to compute the rooted $2$-transitive number of $y$ in $T^{c_i}$. We use another function called \textsc{Mark\_Required$()$}. This function takes the $2$-transitive number of a vertex $x$ and the rooted $2$-transitive number of its children in $T^x$ as input and marks the status of whether a child, say $v$, is required or not to achieve the $2$-transitive number of $x$. We mark $R(v)=1$ if the child $v$ is required; otherwise, $R(v)=0$. We compute the $2$-transitive number of every vertex, other than $c$, by processing the vertices in the reverse order of $\sigma$, that is, in a top-down approach in $T^c$. While processing the vertex $c_i$, first, based on the status marked by the \textsc{Mark\_Required$()$} function, we calculate the rooted transitive number of $p(c_i)$ in $T^{c_i}$, where $p(c_i)$ is the parent of $c_i$ in the rooted tree $T^c$. Then, we call \textsc{$2$-Transitive\_Number$()$} to calculate the $2$-transitive number of $c_i$. Next, we call the \textsc{Mark\_Required$()$} to mark the status of the children, which will be used in subsequent iterations. At the end of this top-down approach, we have the $2$-transitive number of all the vertices, and hence the $2$-transitive number of the tree $T$. The process of finding $Tr_2(T)$ is described in Algorithm \ref{Algo:2-trasitivity(T)}.

\subsubsection{Proof of correctness}

In this subsection, we prove the correctness of Algorithm \ref{Algo:2-trasitivity(T)}. It is clear that the correctness of Algorithm \ref{Algo:2-trasitivity(T)} depends on the correctness of the functions $2$-Transitive\_Number$()$ and Mark\_Required$()$. First, we show the following two lemmas, which prove the correctness of the $2$-Transitive\_Number$()$ function.

\begin{lemma}\label{tree_lemma_2_transitivity}
	Let $v$ be a vertex of $T$, and $t_2(v,T)=t$. Then there exists a $2$-transitive partition of $T$, say $\{V_1, V_2, \ldots, V_i\}$, such that $v\in V_i$, for all $1\leq i\leq t$.
	
\end{lemma}

\begin{proof}
	Since $t_2(v,T)=t$, there is a $2$-transitive partition $\pi=\{U_1, U_2, \ldots, U_t\}$ of $T$ such that $v\in U_t$. For each $1\leq i\leq t$, let us define another $2$-transitive partition $\pi'=\{V_1,V_2,\ldots,V_i\}$ of $T$ as follows: $V_j=U_j$ for all $1\leq j\leq (i-1)$ and $V_i= \displaystyle{\bigcup_{j=i}^{t}U_j}$. Clearly, $\pi'$ is a $2$-transitive partition of $T$ of size $i$ such that $v\in V_{i}$. Hence, the lemma follows.
\end{proof}

%Let $\pi'=\{U_1,U_2,\ldots,U_i, \ldots, U_{t-1}\}$  be  a vertex partition, such that  $U_s=V_s$ for $1\leq s\leq i-1 $ and $U_i=V_i\cup V_{i+1}$ and $U_q=V_{q+1}$ for all $q\geq i+1$. Clearly $\pi'$ is a $2$-transitive partition of order $t-1$ and $v\in U_{t-1}$. Similarly, we construct a $2$-transitive partition $\pi''=\{W_1, W_2, \ldots, W_i, \ldots, W_{t-2}\}$ from $\pi'$ of order $t-2$ and $v$ in $W_{t-2}$. Hence, the lemma follows.

\begin{theorem}\label{tree_theorem_2_transitivity}
	Let $v_1, v_2, \ldots, v_k$ be the children of $x$ in a rooted tree $T^x$, and for each $1\leq i\leq k$, $l_i$ denotes the rooted $2$-transitive number of $v_i$ in $T^x$ with $l_1\leq l_2\leq \ldots\leq l_k$. Let $z$ be the largest integer such that there exists a subsequence of $\{l_i: 1\leq i\leq k\}$, say $(l_{i_1}\leq l_{i_2}\leq \ldots \leq l_{i_{2z}})$ such that $l_{i_{2p-1}}\geq p$ and $l_{i_{2p}}\geq p$, for all $1\leq p\leq z$. Then the $2$-transitive number of $x$ in the underlying tree $T$ is $1+z$, that is, $t_2(x, T)=1+z$.
	
	%
	%
	%
	%
	%		Let $c$ be the root vertex of a tree $T$ and $c$ be adjacent to vertices $v_1, v_2, \ldots, v_k$ and $T_i$ denote the subtrees of $T^c$ rooted at $v_i$, that is $T_i$ is $T_{v_i}^c$. Let $t_2^r(v_i,T^c)=l_i$ and assume $l_1\leq l_2\leq \ldots\leq l_k$. Let $z$ be the largest integer such that there exists a subsequence of $\{l_i: 1\leq i\leq k\}$, say $l_{i_1}\leq l_{i_2}\leq \ldots \leq l_{i_{2z}}$ such that $l_{i_{2p-1}} ~and ~ l_{i_{2p}}\geq p$, for all $1\leq p\leq z$, then $t_2(c, T)=1+z$.  
	
\end{theorem}

\begin{proof}
	For each $1\leq j\leq 2z$, let us consider the subtrees $T_{v_{i_j}}^x$. It is also given that $t_2^r(v_{i_j},T^x)=l_{i_j}$, for $j\in \{1, 2, \ldots, {2z}\}$. For all $1\leq p\leq z$, since $l_{i_{2p-1}},l_{i_{2p}}\geq p$, by Lemma \ref{tree_lemma_2_transitivity}, we know that there exist two $2$-transitive partitions $\pi^{2p-1}=\{V_1^{2p-1}, V_2^{2p-1}, \ldots, V_p^{2p-1}\}$ and $\pi^{2p}=\{V_1^{2p}, V_2^{2p}, \ldots, V_p^{2p}\}$ of $T_{v_{i_{2p-1}}}^x$ and $T_{v_{i_{2p}}}^x$, respectively, such that $v_{i_{2p-1}}\in V_p^{2p-1}$ and $v_{i_{2p}}\in V_p^{2p}$. Let us consider the partition of $\pi=\{V_1, V_2, \ldots, V_z, V_{z+1}\}$ of $T$ as follows: $V_i=\displaystyle{\bigcup_{j=2i-1}^{2z}V_i^j}$, for $2\leq i\leq z$, $V_{z+1}=\{x\}$ and every other vertices of $T$ are put in $V_1$. Clearly, $\pi$ is a $2$-transitive partition of $T$, and therefore $t_2(x,T)\geq 1+z$.
	
	Next, we show that $t_2(x, T)$ cannot be more than $1+z$. If possible, let $t_2(x,T)\geq 2+z$. Then, by Lemma \ref{tree_lemma_2_transitivity}, we have that there exists a $2$-transitive partition $\pi=\{V_1, V_2, \ldots, V_{2+z}\}$ such that $x\in V_{2+z}$. This implies that for each $1\leq i\leq 1+z$, $V_i$ contains two neighbours of $x$, say $v_i$ and $v_{i'}$, such that the rooted $2$-transitive number of both $v_i$ and $v_{i'}$ is greater or equal to $i$, that is, $l_i\geq i$ and $l_{i'}\geq i$. The set $\{l_i, l_{i'}| 1\leq i\leq 1+z\}$ forms a desired subsequence of $\{l_i: 1\leq i\leq k\}$, contradicting the maximality of $z$. Hence, $t_2(x,T)=1+z$.	\end{proof}

%Let us consider the subtrees $T_{v_{i_j}}^x$, for all $v_{i_j}$, where $1\leq j\leq 2z$. So, $t_2^r(v_{i_j},T^x)=l_{i_j}$, for $j\in \{i_1, i_2, \ldots, i_{2z}\}$. Since $l_{i_{2p-1}} ~and ~ l_{i_{2p}}\geq p$, for all $1\leq p\leq z$, then by Lemma \ref{tree_lemma_2_transitivity}, there exits a $2$-transitive partition of $T_{v_{i_{2j-1}}}^x$, say $\{V_1^{2j-1}, V_2^{2j-1}, \ldots, V_j^{2j-1}\}$ for all $1\leq j\leq z$, such that $V_j^{2j-1}=\{v_{i_{2j-1}}\}$. Similarly for  $T_{v_{i_{2j}}}$, we have a $2$-transitive partition $\{V_1^{2j}, V_2^{2j}, \ldots, V_j^{2j}\}$, for all $1\leq j\leq z$, such that $V_j^{2j}=\{v_{i_{2j}}\}$. 
%
%
%Let $\pi=\{V_1, V_2, \ldots, V_z, V_{z+1}\}$ be a vertex partition of $T$, such that $V_i=\displaystyle{\bigcup_{j=2i-1}^{2z}V_i^j}$, for $1\leq i\leq z$ and $V_{z+1}=\{x\}$. Clearly $\pi$ is a $2$-transitive partition of $T$. So, $t_2(x,T)\geq 1+z$. Now, we claim that $t_2(x,T)>1+z$, not possible. If $t_2(x,T)>1+z$, assume $t_2(x,T)=k$, then there exists a $2$-transitive $k$-partition of $T$, say $\pi'=\{V_1, V_2, \ldots, V_k\}$ such that $v\in V_k$. Let assume $v_i, v_i'\in V_i$, $1\leq i\leq k-1$ are children of $x$. Since, $\pi'$ is a $2$-transitive partition and $v_i, v_i'\in V_i$, $1\leq i\leq k-1$, then $t_2^r(v_i,T)\geq i$ and also $t_2^r(v_i',T)\geq i$ . Therefore, we can create a longer subsequence, which contradicts the maximality of $z$. Hence, $t_2(x,T)=1+z$.

Note that in line $6$ of Algorithm \ref{Algo:2-trasitivity(T)}, when $2$-Transitive\_Number$()$ is called, it returns the $2$-transitive number of $c_i$ in $T^c_{c_i}$, which is in fact the rooted $2$-transitive number of $c_i$ in $T^c$. And in line $13$ of Algorithm \ref{Algo:2-trasitivity(T)}, when $2$-Transitive\_Number$()$ is called, then it returns the $2$-transitive number of $c_i$ in $T^{c_i}$, which is the same as $t_2(c_i, T)$. From Lemma \ref{tree_lemma_2_transitivity} and Theorem \ref{tree_theorem_2_transitivity}, we have the following algorithm.

\begin{algorithm}[h]
	\footnotesize
	\caption{\textsc{$2$-Transitive\_Number$(t^r_2(v_{1}, T^x), t^r_2(v_{2}, T^x), \ldots, t^r_2(v_{k}, T^c) )$} } \label{Algo:2transNumber(x)}
	
	\textbf{Input:} Rooted $2$-transitive numbers of children of $x$ with $t^r_2(v_{1}, T^x)\leq t^r_2(v_{2}, T^x)\leq \ldots \leq t^r_2(v_{k}, T^x)$
	
	\textbf{Output:} $2$-transitive number of $x$ in the underlying tree $T$, that is, $t_2(x, T)$.
	
	\begin{algorithmic}[1]

		\State $t_2(x, T)$ $\leftarrow$ 1.
		
		\State $j=0$
		
		\ForAll {$i\leftarrow 1$ to $k$}
		
		\If {$t_2^r(v_i,T^x)\geq t_2(x, T)$}
		
		\State $j=j+1$
		
		\If{$j=2$}
		
		\State $t_2(x, T)=t_2(x, T)+1$.
		
		\State $j=0$
		
		\EndIf
		
		%		   \Else
		%			
		%	        \State  $t_2(x, T)=t_2(x, T)$
		
		\EndIf

		\EndFor
		
		\State \Return($t_2(x, T)$).
		
	\end{algorithmic}
	
\end{algorithm}

Next, we prove the correctness of Mark\_Required$()$. Let $T^x$ be a rooted tree and $t_2(x, T)=z$. A child $v$ of $x$ is said to be required if the $t_2(x, T^x \setminus T_v^x)=z-1$. The function returns the required status of every child of $x$ by marking $R(v)=1$ if required and $R(v)=0$ otherwise. The children of $x$ that are required can be identified using the following theorem.

\begin{theorem} \label{2_transitivity_required_descendants}
	Let $T^x$ be a tree rooted at $x$ and $v_1, v_2, \ldots, v_k$ be its children in $T^x$. Also, let the $2$-transitive number of $x$ be $z$, and for each $1\leq i\leq k$, let $l_i$ denote the rooted $2$-transitive number of $v_i$ in $T^x$ with $l_1\leq l_2\leq \ldots\leq l_k$. Then the following holds:
	
	\begin{enumerate}
		\item[(a)] If $k=2z-2$, then $R(v_i)=1$ for all $1\leq i\leq k$.
		
		\item[(b)] Let $k> 2z-2$. Then for all $1 \leq i \leq k-2z+2$, $R(v_i)=0$.

		\item[(c)] Let $k> 2z-2$. Also, let $k-2z+3\leq i\leq k$. If for all $j$, $k-2z+3\leq j\leq i$, $l_{j-1}\geq \ceil{\frac{j-(k-2z+2)}{2}}$, then $R(v_i)=0$.
		
		\item[(d)] Let $k> 2z-2$ and $k-2z+3\leq i\leq k$. If there exists $j$ in $k-2z+3\leq j\leq i$ such that $l_{j-1}< \ceil{\frac{j-(k-2z+2)}{2}}$, then $R(v_i)=1$.
		
	\end{enumerate}

\end{theorem}

\begin{proof}
	$(a)$ Since $t_2(x, T)=z$, by the Lemma \ref{tree_lemma_2_transitivity}, there exists a $2$-transitive partition of $T$, say $\pi=\{V_1, V_2, \ldots, V_{z}\}$, such that $x\in V_z$. In that case, all the vertices in $\{v_1, v_2, \ldots, v_k\}$ must be in $V_1, V_2, \ldots, V_{z-1}$, and each set $V_i$ contains at least two of these vertices. Since $k=2z-2$, each set $V_i$ contains exactly two vertices from $\{v_1, v_2, \ldots, v_k\}$. Therefore, if we remove any $v_i$ from the tree, the $2$-transitive number of $x$ will decrease by $1$. Hence, every $v_i$ is required, that is, $R(v_i)=1$ for all $1\leq i\leq k$.
	
	%, $v_1, v_2, \ldots, v_k$ must be assigned $1, 1, 2, 2, \ldots, z-1, z-1$, respectively, to achieve the value of the parent. Now, if we remove any $v_i$ from the sequence, then $z$ will be decreased by $1$. Hence, each child is a ``Required'' child. Therefore, $R(v_i)=1$, for all $1\leq i\leq z$.
	$(b)$ Let $\pi=\{V_1, V_2, \ldots, V_{z}\}$ be a $2$-transitive partition of $T$ such that $x\in V_z$. In this partition, at least two vertices from $\{v_1, v_2, \ldots, v_k\}$ must be in each $V_i$ for $1\leq i\leq z-1$. As the vertices are arranged in increasing order of their rooted $2$-transitive number, without loss of generality, we can assume that $\{v_1, v_2, \ldots, v_{k-2z+2}\}\subset V_1$ and $v_i\in V_{I_i}$ for each $k-2z+3\leq i \leq k$, where $I_i=\ceil{\frac{i-(k-2z+2)}{2}}$. If we remove any $v_i$ ($1\leq i\leq k-2z+2$) from the tree, the $2$-transitive number of $x$ will be unchanged. Hence, for each $1\leq i\leq k-2z+2$, $v_i$ is not required, that is, $R(v_i)=0$ for all $1\leq i\leq k-2z+2$.
	
	$(c)$ Let us consider the same $2$-transitive partition $\pi$ of $T$ as in case $(b)$. Let for some $k-2z+3\leq i \leq k$,  $v_i$ be a vertex such that $l_{j-1}\geq I_j$ for all $k-2z+3\leq j\leq i$, where $I_j= \ceil{\frac{j-(k-2z+2)}{2}}$. In this case, we can modify $\pi$ to get a $2$-transitive partition of $T\setminus \{v_i\}$ of size $z$. The modification is as follows: for each $j\in \{k-2z+3, k-2z+4, \ldots, i\}$, we put $v_{i-1}\in V_{I_i}$ and remove the vertices that are not in $T\setminus \{v_i\}$. Therefore, $v_i$ is not required and $R(v_i)=0$ for such vertices.
	
	%%%Assume $k>2z-2$. Then the following assignment of the vertices $v_i$ gives us a $2$-transitive partition of size $z$ for $c$. For $1\leq i\leq k-2z+2$, assign $v_i$ to $1$, and for $k-2z+2<i\leq k$, assign $v_i$ to $j_{i}$, where $j_i=\ceil{\frac{i-(k-2z+2)}{2}}$. Clearly, if we remove $v_i$, $1\leq i\leq k-2z+2$, then there will be no changes in $t_2(c, T^c)$. Therefore, $v_i$ is ``Not Required'', for all $1\leq i\leq k-2z+2$, which implies $R(v_i)=0$. Consider any $i$, $k-2z+1<i\leq k$. If $l_{j-1}\geq \ceil{\frac{j-(k-2z+2)}{2}}$ for all $j$, $k-2z+2<j\leq i$, then $v_i$ can be removed and there will be no changes in $t_2(c,T^c)$, which is obtained by assigning the vertices $\{v_{l}, v_{l+1}, \ldots, v_{i-1}, v_{i+1}, \ldots, v_k\}$, $1, 1, 2, 2, \ldots, z-1, z-1$, respectively. So, $v_i$ is ``Not Required'' which implies $R(v_i)=0$.
	$(d)$ Let for some $i$, $k-2z+3\leq i \leq k$, and $v_i$ be a vertex such that $l_{j-1}< I_j$ for some $k-2z+3\leq j\leq i$, where $I_j= \ceil{\frac{j-(k-2z+2)}{2}}$. Since the rooted $2$-transitive numbers are arranged in increasing order, $l_p< I_j$ for all $1\leq p\leq j-1$. Suppose, after deleting the vertex $v_i$, the $2$-transitive number of $x$ in $T\setminus \{v_i\}$ remains $z$. Let $\pi'=\{V_1, V_2, \ldots, V_{z}\}$ be a $2$-transitive partition of $T\setminus \{v_i\}$ of size $z$ such that $x\in V_{z}$. Since $l_p< I_j$ for all $1\leq p\leq j-1$, none of the vertices of $\{v_1, v_2, \ldots, v_{j-1}\}$ can be in the sets $V_{I_j}, V_{I_{j+1}}, \ldots, V_{I_k}$. On the other hand, the sets $V_{I_j}, V_{I_{j+1}}, \ldots, V_{I_k}$ must contain at least $(k-j+1)$ vertices from $\{v_1, v_2, \ldots, v_{i-1}, v_{i+1}, \ldots, v_{k}\}$, as $\pi'$ is a $2$-transitive partition of $T\setminus \{v_i\}$. Therefore, $V_{I_j}, V_{I_{j+1}}, \ldots, V_{I_k}$ contains at least $(k-j+1)$ vertices from $\{v_{j}, v_{j+1}, \ldots, v_{i-1}, v_{i+1}, \ldots, v_{k}\}$. But there are only $(k-j)$ many vertices available. Hence, the $2$-transitive number of $x$ in $T\setminus \{v_i\}$ cannot be $z$. Therefore, $v_i$ is required and $R(v_i)=1$ for such vertices.
\end{proof}

Note that the condition in case $(d)$ is such that if $R(v_i)=1$ for some $i$, then $R(v_j)=1$ for all $i+1\leq j\leq k$. Based on the Theorem \ref{2_transitivity_required_descendants}, we have the following function.

\begin{algorithm}[h]
	\footnotesize
	
	\caption{\textsc{Mark\_Required$(t_2(x, T), t^r_2(v_1, T^x), t^r_2(v_2, T^x), \ldots, t^r_2(v_k, T^x) )$}}  \label{Algo:2required_descendants(c)}
	
	\textbf{Input:} A rooted tree $T^c$, rooted at a vertex $c$ and rooted $2$-transitive number of the children of $c$, such that $t^r_2(v_1, T^c)\leq t^r_2(v_2, T^c)\leq \ldots \leq t^r_2(v_k, T^c)$.
	
	\textbf{Output:} $R(v)$ value of $v$, $v$ is a child of $x$ in $T^x$.
	
	\begin{algorithmic}[1]
		
		\If {the number of children is $2z-2$, that is $k=2z-2$}
		
		\State $R(v)=1$, for all children $v$ of $x$ in $T^x$.
		
		\EndIf
		\If {$k> 2z-2$}
		
		\State For all $1\leq i\leq k-2z+2$, $R(v_i)=0$.
		
		\EndIf
		
		\ForAll {$i\leftarrow k-2z+3$ to $k$}
		
		\If {$l_{i-1}\geq \ceil{\frac{i-(k-2z+2)}{2}}$}
		
		\State $R(v_i)=0$
		
		\Else 
		
		\State  $R(v_i)=1$
		
		\State {\bf break}
		
		\EndIf

		\EndFor

		\State For all $j> i$, $R(v_j)=1$
		
		%	\ForAll {$k-2z+2<j\leq i$ }
		%		
		%		\State $l_{j-1}\geq \ceil{\frac{j-(k-2z+2)}{2}}$
		%		
		%		\State $R(v_i)=0$.
		%		
		%		\EndFor
		%       \EndFor
		%		\State $R(v_i)=1$
		
	\end{algorithmic}
	
\end{algorithm}

\subsubsection{Complexity Analysis}
In the function \textsc{$2$-Transitive\_Number()}, we find the $2$-transitive number of vertex $x$ based on the rooted $2$-transitive number of its children. We assume the children are sorted according to their rooted $2$-transitive number. Since the for loop in lines $3-8$ of \textsc{$2$-Transitive\_Number()} runs for every child of $x$, this function takes $O(deg(x))$ time. Similarly, the function \textsc{Mark\_Required()} takes the $2$-transitive number of a vertex $x$ and the rooted $2$-transitive number of its children in $T^x$ as input and marks the status of whether a child, say $v$, is required or not to achieve the $2$-transitive number of $x$. Here also, we have assumed that the children are sorted according to their rooted $2$-transitive number. Clearly, lines $3-4$ of \textsc{Mark\_Required()} take $O(deg(x))$. In lines $3-4$, we mark the status for a few children without any checking, and for each of the remaining vertices, we mark the required status by checking the condition in $O(1)$ time. Therefore, \textsc{Mark\_Required()} also takes $O(deg(x))$ time. In the main algorithm \textsc{$2$-Transitivity(T)}, the vertex ordering mentioned in line $1$ can be found in linear time. Then, in a bottom-up approach, we calculate the rooted $2$-transitive numbers of every vertex. For that, we are spending $O(deg(c_i))$ for every $c_i\in \sigma$. Note that we must pass the children of $c_i$ in a sorted order to \textsc{$2$-Transitive\_Number()}. But as discussed in \cite{hedetniemi1982linear} (an algorithm for finding the Grundy number of a tree), we do not need to sort all the children based on their rooted $2$-transitive numbers; sorting the children whose rooted $2$-transitive number is less than $deg(c_i)$ is sufficient. We can argue that this can be done in $O(deg(c_i))$, as shown in \cite{hedetniemi1982linear}. Hence, the loop in lines $1-6$ takes linear time. In similar ways, we conclude that lines $8-14$ take linear time. Therefore, we have the following theorem:

\begin{theorem}
	The \textsc{M$2$TP} can be solved in linear time for trees.
\end{theorem}

% and \textsc{Mark\_Required}, for a vertex $v$, we process a list of children of length at  most $deg(v)$. The condition in the algorithm is automatically satisfied if the rooted $2$-transitive number of a child is more or equal to $deg(v)$. If we assume the list is already sorted, then the total amount of time required to process this list in each procedure is at most $deg(v)$. Now the algorithm \textsc{$2$-Transitivity $Tr_2(T)$}  can be computed in $O(n)$ time. As the for loop in line $4-8$ runs in $O(n)$, because $\displaystyle{\sum_{v\in T}deg(v)}=2(n-1)$. Also, the for loop in line $10-14$ runs in $O(n)$, because $\displaystyle{\sum_{v\in T}deg(v)}=2(n-1)$ and the line can be computed in $O(n)$. From the analysis given in  \cite{hedetniemi1982linear}, we are required to sort the list up to $deg(v)$; hence, sorting will take at most $O(deg(v))$ time for a vertex $v$.

\subsection{Split graphs}
In this subsection, we propose a linear-time algorithm for finding the $2$-transitivity of a given split graph. We first show that the $2$-transitivity of  a split graph $G$ can be either $\ceil{\frac{\omega(G)}{2}}$ or $\ceil{\frac{\omega(G)}{2}}+1$, where $\omega(G)$ is the size of a maximum clique in $G$. Recall that a graph $G=(V, E)$ is a split graph if $V$ can be partitioned into an independent set $S$ and a clique $K$.

\begin{lemma}\label{split_graph_2_transitivity_lemma_1}
	Let $G=(S\cup K, E)$ be a split graph, where $S$ and $K$ are the independent set and clique of $G$, respectively, such that $\omega(G)=|K|$. Then $\ceil{\frac{\omega(G)}{2}}\leq Tr_2(G)\leq \ceil{\frac{\omega(G)}{2}}+1$.
\end{lemma}

\begin{proof}
	Consider a vertex partition, $\pi=\{V_1, V_2, \ldots, V_{\ceil{\frac{\omega(G)}{2}}}\}$ of $G$ as follows: if $\omega(G)$ is even, then each $V_i$ contains exactly two vertices from $K$ for $1\leq i\leq \ceil{\frac{\omega(G)}{2}}$, and every other vertex of $G$ is in $V_1$. If $\omega(G)$ is odd, then each $V_i$ contains exactly two vertices from $K$ for $1\leq i\leq \ceil{\frac{\omega(G)}{2}}-1$, and $V_{\ceil{\frac{\omega(G)}{2}}}$ contains exactly one vertex of $K$ and every other vertex of $G$ is in $V_1$. It can be easily verified that $\pi$ forms a $2$-transitive partition of size $\ceil{\frac{\omega(G)}{2}}$. Hence, $\ceil{\frac{\omega(G)}{2}}\leq Tr_2(G)$.
	
	We prove the upper bound by contradiction. Let $Tr_2(G)\geq \ceil{\frac{\omega(G)}{2}}+2$ and $\pi=\{V_1, V_2, \ldots, V_{\ceil{\frac{\omega(G)}{2}} +2}\}$ be a $2$-transitive partition of $G$ with size $\ceil{\frac{\omega(G)}{2}}+2$. First, note that there do not exist two sets in $\pi$ that contain vertices from $S$ only. Therefore, all sets, except one, of $\pi$ contain vertices from $K$. Now, since $|K|=\omega(G)$, there exist at least two sets in $\pi$, say $V_i$ and $V_j$ ($i<j$), such that $V_i$ and $V_j$ contain at least one vertex from $S$ but at most one vertex from $K$. In that case, $V_j$ is not $2$-dominated by $V_i$ as $S$ is an independent set. Hence, $Tr_2(G)\leq \ceil{\frac{\omega(G)}{2}}+1$. Therefore, $\ceil{\frac{\omega(G)}{2}}\leq Tr_2(G)\leq \ceil{\frac{\omega(G)}{2}}+1$. \end{proof}

%
%Now consider a vertex partition, $\pi=\{V_1, V_2, \ldots, V_{\ceil{\frac{\omega(G)}{2}}}\}$ of $G$. If $\omega(G)$ is even then consider each $V_i$, $1\leq i\leq \ceil{\frac{\omega(G)}{2}}$ contains exactly two vertices from $K$ and all the other vertex of $G$ in $V_1$ and if $\omega(G)$ is odd then consider each $V_i$, $1\leq i\leq \ceil{\frac{\omega(G)}{2}}-1$ contains exactly two vertices from $K$ and $V_{\ceil{\frac{\omega(G)}{2}}}$ contains exactly one vertex of $K$ and all the other vertex of $G$ in $V_1$. Clearly $\pi$ is  a $2$-transitive partition of size  $\ceil{\frac{\omega(G)}{2}}$. So, $Tr_2(G)\geq \ceil{\frac{\omega(G)}{2}}$. Hence, $\ceil{\frac{\omega(G)}{2}}\leq Tr_2(G)\leq \ceil{\frac{\omega(G)}{2}}+1$.		 Since $|K|=\omega(G)$, then there exist at least two partitions in $\pi$, say $V_i$ and $V_j$ with $i<j$, such that $V_i$ and $V_j$ contains at least one vertex from $S$ and at most one vertex from $K$. In this case, $V_j$ is not $2$-dominated by $V_i$ as $S$ is an independent set. Therefore, $Tr_2(G)\leq \ceil{\frac{\omega(G)}{2}}+1$. 

Next, we characterize the split graphs with $2$-transitive number equal to $\ceil{\frac{\omega(G)}{2}}+1$. The proof is divided into two cases based on the parity of $\omega(G)$.

%The following Lemma \ref{Lemma_split_odd_clique_number} characterize split graphs with $2$-transitive number equal to $\ceil{\frac{\omega(G)}{2}}+1$, when $\omega(G)$ odd, and Lemma \ref{Lemma_split_even_clique_number} characterize split graphs with $2$-transitive number equal to $\frac{\omega(G)}{2}+1$, when $\omega(G)$ even.

\begin{lemma}\label{Lemma_split_odd_even_clique_number}
	Let $G=(S\cup K, E)$ be a split graph, where $S$ and $K$ are the independent set and clique of $G$, respectively, such that $\omega(G)=|K|$.
	
	(a) If $\omega(G)$ is odd, then $Tr_2(G)=\ceil{\frac{\omega(G)}{2}}+1$ if and only if every vertex of $K$ has at least two neighbours in $S$.
	
	(b) If $\omega(G)$ is even, then $Tr_2(G)=\ceil{\frac{\omega(G)}{2}}+1$ if and only if at least $|K|-1$ vertices in $K$ have a neighbour in $S$.
\end{lemma}

\begin{proof}
	(a) Let $Tr_2(G)=\ceil{\frac{\omega(G)}{2}}+1$ and $\pi=\{V_1,V_2, \ldots, V_l, \ldots, V_{\ceil{\frac{\omega(G)}{2}}+1}\}$ be a $2$-transitive partition of $G$. First note that if each set in $\pi$ contains vertices from $K$, then there exist at least two sets in $\pi$, say $V_i$ and $V_j$ $(i< j)$ such that both $V_i$ and $V_j$ contain at least one vertex from $S$ but at most one vertex from $K$, as the size of $K$ is $\omega(G)$. In that case, $V_j$ is not $2$-dominated by $V_i$ as $S$ is an independent set. So, at least one set in $\pi$ exists that contains only vertices from $S$. Note that there do not exist two sets in $\pi$ that contain vertices from $S$ only. Therefore, exactly one set in $\pi$, say $V_l$, contains only vertices from $S$. First, we show that $l$ cannot be $\ceil{\frac{\omega(G)}{2}}+1$. If $l=\ceil{\frac{\omega(G)}{2}}+1$, then each of $V_i$, $1\leq i\leq \ceil{\frac{\omega(G)}{2}}$, contains at least two vertices from $K$. That is not possible as the size of $K$ is $\omega(G)$ and $\omega(G)$ is odd. So, $V_l$ is not the last set in $\pi$, and hence, $|V_l|\geq 2$. As $V_l$ $2$-dominates $V_j$ for all $j>l$, every vertex in $V_j$ has two neighbours in $V_l\subseteq S$. Clearly, $V_j$ does not contain any vertices from $S$, and all the vertices of $K$ that appear in a set $V_j$ ($j>l$) have at least two neighbours in $S$. Also, note that each set $V_i$ $(i< l)$ contains exactly two vertices from $K$. Otherwise, we have the previous case, where two sets contain at least one vertex from $S$ but at most one vertex from $K$. Since $V_i$ $2$-dominates $V_l$ and contains exactly two vertices from $K$, each vertex from $K \cap V_i$ has two neighbours in $V_l\subseteq S$. Hence, all the vertices of $K$ are adjacent to at least two vertices of $S$.

	Conversely, let $G=(S\cup K, E)$ be a split graph that satisfies the condition; every vertex of $K$ is adjacent to at least two vertices of $S$. Now consider a vertex partition, $\pi=\{V_1, V_2, \ldots, V_{\ceil{\frac{\omega(G)}{2}}+1}\}$ of $G$ such that $V_1=S$, each $V_i$ for $2\leq i\leq \ceil{\frac{\omega(G)}{2}}$ contains exactly two vertices from $K$ and $V_{\ceil{\frac{\omega(G)}{2}}+1}$ contains exactly one vertex from $K$. Since all the vertices of $K$ are adjacent to at least two vertices of $S$, $V_1$ $2$-dominates $V_j$ for all $j\geq 2$. Also, it is easy to verify that $V_i$ $2$-dominates $V_j$ for all $2\leq i<j \leq \ceil{\frac{\omega(G)}{2}}+1$ as $K$ is a clique. Therefore, $\pi$ is a $2$-transitive partition of size $\ceil{\frac{\omega(G)}{2}}+1$. By the Lemma \ref{split_graph_2_transitivity_lemma_1}, it follows that $Tr_2(G)=\ceil{\frac{\omega(G)}{2}}+1$.

	% So, $\pi$ is $2$-transitive partition of size $\ceil{\frac{\omega(G)}{2}}+1$. Therefore, $Tr_2(G)\geq \ceil{\frac{\omega(G)}{2}}+1$. Again by the Lemma \ref{split_graph_2_transitivity_lemma_1}, $Tr_2(G)\leq  \ceil{\frac{\omega(G)}{2}}+1$. So, $Tr_2(G)=\ceil{\frac{\omega(G)}{2}}+1$.

	(b) Let $Tr_2(G)=\ceil{\frac{\omega(G)}{2}}+1$ and $\pi=\{V_1,V_2, \ldots, V_l, \ldots, V_{\ceil{\frac{\omega(G)}{2}}+1}\}$ be a $2$-transitive partition of $G$. If every set in $\pi$ contains vertices from $K$, then by the pigeonhole principle, at least two sets contain exactly one vertex from $K$, as $|K|=\omega(G)$. Clearly, one of these two sets must contain at least one vertex from $S$. On the other hand, if every set in $\pi$ does not contain vertices from $K$, then exactly one set contains vertices from $S$ only as $S$ is an independent set. We prove the lemma by dividing the proof into the cases mentioned above.
	
	\begin{case}
		There is a set, say $V_p$, in $\pi$ containing at least one vertex from $S$ and exactly one vertex from $K$ 
	\end{case}
	In this case, note that every set $V_j$ ($j>p$) contains no vertex from $S$. Since $V_p$ $2$-dominates $V_j$ $(j>p)$, each vertex from $V_j$ must have a neighbour in $S$. Also, to $2$-dominate $V_p$, each set $V_i$ $(i< p)$ contains at least two vertices from $K$. Since $|K|=\omega(G)$, it follows that each $V_i$ $(i< p)$ must contain exactly two vertices from $K$. Therefore, at least $|K|-1$ vertices in $K$ have a neighbour in $S$.
	
	\begin{case}
		There exists exactly one set, say $V_l$, in $\pi$ containing only vertices from $S$
	\end{case}
	
	In this case, note that every set $V_j$ ($j>l$) contains no vertex from $S$. So, each $V_j$ for $l+1 \leq j\leq \ceil{\frac{\omega(G)}{2}}$ contains at least two vertices from $K$, and $V_{\ceil{\frac{\omega(G)}{2}}+1}$ contains at least one vertex from $K$. Also, to $2$-dominate $V_l$, each $V_i$ for $1\leq i\leq l-1$ contains at least two vertices from $K$. Therefore, the distribution of the vertices of $K$ in $\pi$ is as follows: either each set in $\pi$ contains exactly two vertices from $K$ or one set contains exactly $1$ vertex from $K$, another set contains exactly $3$ vertices from $K$, and the rest contain exactly $2$ vertices from $K$. If each set in $\pi$ contains exactly two vertices from $K$, then every vertex of $K$ has a neighbour in $S$. Otherwise, let $V_t$ be the set that contains exactly $3$ vertices from $K$. If $t>l$, every vertex of $K$ has a neighbour in $S$. Otherwise, at least $|K|-1$ vertices have a neighbour in $S$.
	
	Conversely, let $G=(S\cup K, E)$ be a split graph that satisfies the condition, that is, at least $|K|-1$ vertices in $K$ have a neighbour in $S$. Let $K_1$ be the subset of $K$ of $|K|-1$ vertices with a neighbour in $S$. Now consider a vertex partition $\pi=\{V_1, V_2, \ldots, V_{\ceil{\frac{\omega(G)}{2}}+1}\}$ of $G$, such that $V_1= S \cup (K\setminus K_1)$, each $V_i$ for $2< i\leq \ceil{\frac{\omega(G)}{2}}$ contains exactly two vertices from $K_1$ and $V_{\ceil{\frac{\omega(G)}{2}}+1}$ contains exactly one vertices from $K_1$. Since every vertex of $K_1$ is adjacent to at least one vertex of $S$ and $K$ is a clique, $V_1$ $2$-dominates $V_i$ for all $i\geq 2$. Also, it is easy to verify that $V_i$ $2$-dominates $V_j$ for all $2\leq i<j \leq \ceil{\frac{\omega(G)}{2}}+1$ as $K$ is a clique. Therefore, $\pi$ is a $2$-transitive partition of size $\ceil{\frac{\omega(G)}{2}}+1$. By the Lemma \ref{split_graph_2_transitivity_lemma_1}, it follows that $Tr_2(G)=\ceil{\frac{\omega(G)}{2}}+1$.
\end{proof}

	Based on Lemma \ref{Lemma_split_odd_even_clique_number}, we can design an algorithm for finding $2$-transitivity of a given split graph. As	the conditions mentioned in Lemma \ref{Lemma_split_odd_even_clique_number} can be checked from the degrees of the vertices in $K$, the proposed algorithm runs in linear time.	Hence, we have the following theorem:

	\begin{theorem}
		The \textsc{Maximum $2$-Transitivity Problem} can be solved in linear time for split graphs.
	\end{theorem}

	\subsection{Bipartite chain graphs}
	
	In this subsection, we find the $2$-transitivity of a bipartite chain graph $G$ in terms of its transitivity. A bipartite graph $G=(X\cup Y,E)$ is called a \textit{bipartite chain graph} if there exists an ordering of vertices of $X$ and $Y$, say $\sigma_X= (x_1,x_2, \ldots ,x_m)$ and $\sigma_Y=(y_1, y_2, \ldots, y_n)$, such that $N(x_m)\subseteq N(x_{m-1})\subseteq \ldots \subseteq N(x_2)\subseteq N(x_1)$ and $N(y_n)\subseteq N(y_{n-1})\subseteq \ldots \subseteq N(y_2)\subseteq N(y_1)$. This ordering of $X$ and $Y$ is called a \emph{chain ordering}. A chain ordering of a bipartite chain graph can be computed in linear time \cite{heggernes2007linear}. We first show the following result to find the $2$-transitivity of a bipartite chain graph $G$.

	%design the algorithm, we first prove that if $t$ is the maximum integer such that $G$ contains $K_{2\floor{\frac{t}{2}}, 2\floor{\frac{t}{2}}-1}-\{e, e'\}$ as a subgraph, where $e, e'$ are two independent edges of $G$, then $Tr_2(G)=\floor{\frac{t}{2}}+1$.
	
	%After that, we design an algorithm for finding maximum integer $t$ such that $G$ contains either $K_{t,t}$ or  $K_{t,t}-\{e\}$ as an induced subgraph.

	\begin{lemma} \label{two_transitivity_chain_lemma_1} 
		The $2$-transitivity of $G=(K_{2\floor{\frac{t}{2}}, 2\floor{\frac{t}{2}}-1}\setminus C_4)$ is $\floor{\frac{t}{2}}+1$.
	\end{lemma}
	\begin{proof}
		Let $t$ be an even integer and $V(K_{t, t-1})=X\cup Y$, where the set $X=\{x_1, x_2, \ldots, x_t\}$ and $Y=\{y_1, y_2, \ldots, y_{t-1}\}$. Also, without loss of generality, let us assume that $\{x_t,x_{t-1}, y_{t-1}, y_{t-2}\}$ induces a $C_4$ in $G$. Now consider a vertex partition $\pi=\{V_1, V_2, \ldots, V_{k}\}$ of size $k={\frac{t}{2}}+1$ of $G$ as follows: $V_j=\{x_{t-2j+2}, x_{t-2j+1}, y_{t-2j+1}, y_{t-2j}\}$, for all $1\leq j\leq k-2$, $V_{k-1}=\{x_1, x_2\}$ and $V_k=\{y_1\}$. Clearly, $\pi$ is a $2$-transitive partition of $G$. Therefore, $Tr_2(G)\geq k=\frac{t}{2}+1$. Also, by Proposition \ref{upper_bound_2-transitivity}, we know $Tr_2(G)\leq \floor{\frac{\Delta(G)}{2}}+1=\floor{\frac{t}{2}}+1$. Hence, the $2$-transitivity of $G$ is $\frac{t}{2}+1$, where $t$ is even. Now, let us assume that $t$ is an odd integer and $t'=t-1$. Clearly, the graph $(K_{2\floor{\frac{t}{2}}, 2\floor{\frac{t}{2}}-1}-\setminus C_4)$ is the same as $(K_{t', t'-1}\setminus C_4)$. As $t'$ is even, we have $Tr_2(G)=Tr_2(K_{t', t'-1}\setminus C_4)=\frac{t'}{2}+1=\floor{\frac{t}{2}}+1$. Hence, the $2$-transitivity of $G$ is $\floor{\frac{t}{2}}+1$.	\end{proof}
	
	% consider the graph $K_{2\floor{\frac{t}{2}}, 2\floor{\frac{t}{2}}-1}-\{e, e'\}$. Let $t'=t-1$, clearly the graph $K_{2\floor{\frac{t}{2}}, 2\floor{\frac{t}{2}}-1}-\{e, e'\}$ is same as $K_{t', t'-1}-\{e, e'\}$, where $t'$ is an even integer. From the above arguments, we have $Tr_2(K_{2\floor{\frac{t}{2}}, 2\floor{\frac{t}{2}}-1}-\{e, e'\})=Tr_2(K_{t', t-1}-\{e, e'\})=\frac{t'}{2}+1=\floor{\frac{t}{2}}+1$. Hence, $2$-transitivity of $K_{2\floor{\frac{t}{2}}, 2\floor{\frac{t}{2}}-1}-\{e, e'\}$, where $e, e'$ are two independent edges is $\floor{\frac{t}{2}}+1$. 
	
	Using the above result, we find the $2$-transitivity of a bipartite chain graph in the following theorem. 
	
	\begin{lemma}\label{two_transitivity_chain_lemma_2} 
		Let $G=(X\cup Y, E)$ be a bipartite chain graph and $t$ be the maximum integer such that $G$ contains $(K_{2\floor{\frac{t}{2}}, 2\floor{\frac{t}{2}}-1}\setminus C_4)$ as a subgraph. Then $Tr_2(G)=\floor{\frac{t}{2}}+1$. 
	\end{lemma}
	
	\begin{proof}
		Since $t$ is the maximum integer such that $G=(X\cup Y,E)$ contains $(K_{2\floor{\frac{t}{2}}, 2\floor{\frac{t}{2}}-1}\setminus C_4)$ as a subgraph,  Lemma \ref{two_transitivity_chain_lemma_1} implies that $Tr_2(G)\geq \floor{\frac{t}{2}}+1$. To show that $Tr_2(G)$ cannot be more than $\floor{\frac{t}{2}}+1$, we first show that if $Tr_2(G)=\floor{\frac{m}{2}}+1$, then $G$ contains $(K_{2\floor{\frac{m}{2}}, 2\floor{\frac{m}{2}}-1}\setminus C_4)$ as a subgraph. We prove this by induction on $m$. The base case for $m=2$ follows from Proposition \ref{2_transitivity_geq_2}, as $P_3$ can be viewed as $(K_{2,1}\setminus C_4)$. By induction hypothesis, let us assume that the above statement is true for any graph $G$ with $2$-transitivity is less than $\floor{\frac{m}{2}}+1$. Let $G$ be a bipartite chain graph with $Tr_2(G)=\floor{\frac{m}{2}}+1$ and $\{V_1,V_2, \ldots ,V_{\floor{\frac{m}{2}}+1}\}$ be a $2$-transitive partition of $G$ of size $(\floor{\frac{m}{2}}+1)$. Let $G'=G[V_2\cup V_3\cup \ldots \cup V_{\floor{\frac{m}{2}}+1} ]$. Clearly,  $Tr_2(G')= \floor{\frac{m}{2}}=\floor{\frac{m-2}{2}}+1$. Let $\alpha=\floor{\frac{m-2}{2}}$. By induction hypothesis $G'$ contains $(K_{2\alpha, 2\alpha-1}\setminus C_4)$ as a subgraph. Let $X'=\{x_1', x_2', \ldots , x_{2\alpha}'\}$ and $Y'=\{y_1', y_2', \ldots , y_{2\alpha-1}'\}$ be two sets such that $G'[X'\cup Y']$ also contains $(K_{2\alpha, 2\alpha-1}\setminus C_4)$ as a subgraph. Let $x_{i}', x_{j}'\in X'$ with $i<j$ and $y_{r}', y_{s}'\in Y'$ with $r<s$ be the vertices such that if we remove all possible edges among them from $G'[X'\cup Y']$, then we get a copy of $(K_{2\alpha, 2\alpha-1}\setminus C_4)$. Now, in $G$, since $V_1$ $2$-dominates $V_j$, for all  $j\geq 2$, $V_1$ contains at least two vertices from $X$ and at least two vertices from $Y$. Let $\{x_{l_1},x_{l_2}, \ldots, x_{l_s}\}$ and $\{y_{k_1},y_{k_2}, \ldots, y_{k_t}\}$ be the vertices in $V_1$ from $X$ and $Y$, respectively. Since $G$ is a bipartite chain graph, there exist $x_p, x_{p'}\in\{x_{l_1},x_{l_2}, \ldots, x_{l_s}\}$ and $ y_q, y_{q'}\in \{y_{k_1},y_{k_2},\ldots, y_{k_t}\}$ such that $N(x_p)\supseteq N(x_{p'})\supseteq N(x)$ for all $x\in\{x_{l_1}, \ldots, x_{l_s}\}\setminus \{x_p, x_{p'}\}$ and $N(y_q)\supseteq N(y_{q'})\supseteq N(y)$ for all 
		$y \in\{y_{k_1},y_{k_2}, \ldots, y_{k_t}\}\setminus \{y_q, y_{q'}\}$. Therefore, both $x_p, x_{p'}$ $2$-dominate $\{y_1', y_2', \ldots, y_{2\alpha-1}'\}$ and both $y_q, y_{q'}$ $2$-dominate $\{x_1', x_2', \ldots, x_{2\alpha}'\}$. 
		
		Now we focus on the edges between $\{x_i', x_j'\}$ and $\{y_r', y_s'\}$ and also the edges between $\{x_p, x_{p'}\}$ and $\{y_q, y_{q'}\}$. Since $G$ is a bipartite chain graph, either $N(x_p)\subseteq N(x_i')$ or $N(x_i')\subseteq N(x_p)$. We show the existence of the required subgraph for one case, and the other follows similarly. Let us assume that $N(x_i')\subseteq N(x_p)$. Since $\{y_q, y_{q'}\} \subseteq N(x_i')$, we have $x_py_q , x_py_{q'}\in E$ in this case. Now, if $N(x_i')\subseteq N(x_{p'})$, then we can show similarly that $x_{p'}y_q , x_{p'}y_{q'}\in E$. Therefore, $G[X' \cup Y' \cup \{x_p, x_p',y_q, y_q' \}]$ contains a $(K_{2\alpha+2, 2\alpha+1}\setminus C_4)$ as a subgraph, which is obtained by deleting edges between $\{x_i', x_j'\}$ and $\{y_r', y_s'\}$. Otherwise, let $N(x_{p'})\subseteq N(x_i')$. Since $\{y_r', y_s'\}\subseteq N(x_{p'})$, we have $x_i'y_r', x_i'y_s'\in E$. Further, we consider the relation between $N(x_j')$ and $N(x_p')$. If $N(x_j')\subseteq N(x_{p'})$, then we have $x_{p'}y_q, x_{p'}y_{q'}\in E$ as $\{y_q, y_{q'}\}\subseteq N(x_{j}')$. Therefore, $G[X' \cup Y' \cup \{x_p, x_p',y_q, y_q' \}]$ contains a $(K_{2\alpha+2, 2\alpha+1}\setminus C_4)$ as a subgraph, which is obtained by deleting edges between $\{x_i', x_j'\}$ and $\{y_r', y_s'\}$. On the other hand, if $N(x_{p'})\subseteq N(x_j')$, then we have $x_{j}'y_r', x_{j}'y_{s}'\in E$ as $\{y_r', y_{s}'\}\subseteq N(x_{p'})$. Therefore, $G[X' \cup Y' \cup \{x_p, x_p',y_q, y_q' \}]$ contains a $(K_{2\alpha+2, 2\alpha+1}\setminus C_4)$ as a subgraph, which is obtained by deleting edges between $\{x_p, x_{p'}\}$ and $\{y_{q}, y_{q'}\}$. Hence, for every case, we show that if $Tr_2(G)=\floor{\frac{m}{2}}+1$, then $G=(X\cup Y,E)$ contains $(K_{2\floor{\frac{m}{2}},2\floor{\frac{m}{2}}-1}\setminus C_4)$ as a subgraph. If $Tr_2(G)\geq \floor{\frac{t}{2}}+2$, then $G$ contains $(K_{2\floor{\frac{t}{2}}+1,2\floor{\frac{t}{2}}}\setminus C_4)$ as a subgraph, which contradicts the maximality of $t$. Therefore, $Tr_2(G)= \floor{\frac{t}{2}}+1$.		\end{proof}

	The above lemma shows that for computing the $2$-transitivity of a bipartite chain graph $G$, we have to check the existence of $(K_{2\floor{\frac{t}{2}}, 2\floor{\frac{t}{2}}-1}\setminus C_4)$ as a subgraph in $G$ for the maximum value of such $t$. It is known that transitivity in a bipartite chain graph $G$ is $t+1$ if and only if $G$ contains either $K_{t,t}$ or $K_{t,t}\setminus\{e\}$ as an induced subgraph for the maximum value of such $t$ \cite{paul2023transitivity}. By comparing the existence of such subgraphs, we show the relation between $Tr_2(G)$ and $Tr(G)$.

	\begin{theorem}
		Let $G=(X\cup Y,E)$ be a bipartite chain graph, and  $\sigma_X= (x_1,x_2, \ldots ,x_m)$ and $\sigma_Y=(y_1,y_2, \ldots ,y_n)$ be the chain ordering of $G$. Also, assume that $t$ is the maximum integer such that $G$ contains either $K_{t,t}$ or $K_{t,t}\setminus\{e\}$ as an induced subgraph. Then the $2$-transitivity of $G$ is given as follows:
		
		\[
		Tr_2(G) = 
		\begin{cases}
			\floor{\frac{Tr(G)}{2}}+1 & if ~G \text{ contains } K_{t, t} \text{ where $t$ is maximum such }\\ & \text{integer and } x_ty_{t+1}\in E \text{ or } x_{t+1}y_t\in E \text{ but not }\\ & \text {both} \\\\ 
			
			\floor{\frac{Tr(G)-1}{2}}+1 & if ~G \text{ contains } K_{t, t}  \text{  where $t$ is maximum such  }\\ & \text {integer and }x_ty_{t+1}\notin E \text{ and } x_{t+1}y_t\notin E \\\\
			
			\floor{\frac{Tr(G)-1}{2}}+1 & if ~G \text{ contains } K_{t, t}\setminus \{e\} \text{  where $t$ is maximum }\\ & \text{ such integer}
			
			%	\floor{\frac{Tr(G)}{2}}+1 &  x_ty_{t+1}\in E \text{ and } x_{t+1}y_t\in E 
			
		\end{cases}
		\]

	\end{theorem}
	
	\begin{proof}
		Let $t$ be the maximum integer such that $G$ contains $K_{t,t}$ as an induced subgraph. From \cite{paul2023transitivity}, we have $Tr(G)=t+1$ and $K_{t, t}$ is induced by $\{x_1, x_2, \ldots, x_t\}\cup \{y_1, y_2, \ldots, y_t\}$. Without loss of generality, assume $x_{t+1}y_{t}\in E$ and $x_{t}y_{t+1}\notin E$. Therefore, $G$ contains $K_{t+1, t}$ as a subgraph. Again, $x_{t+j}y_{t+1}\notin E$, for $j\geq 2$ as $G$ is a bipartite chain graph and $N(x_t)$ does not contain $y_{t+1}$. Also, from our assumption, it follows that $x_{t+1}y_{t+1}\notin E$. So, in this case, $G$ contains $(K_{2\floor{\frac{t'}{2}}, 2\floor{\frac{t'}{2}}-1}\setminus C_4)$ as a subgraph for maximum $t'$, where $t'=t+1$, if $t$ is even and $t'=t+2$ if $t$ is odd. Now, by the Lemma \ref{two_transitivity_chain_lemma_2}, we have $Tr_2(G)=\floor{\frac{t'}{2}}+1=\floor{\frac{t+1}{2}}+1=\floor{\frac{Tr(G)}{2}}+1$, when $t$ is even. Otherwise, $Tr_2(G)=\floor{\frac{t'}{2}}+1=\floor{\frac{t+2}{2}}+1=\frac{t+1}{2}+1=\frac{Tr(G)}{2}+1=\floor{\frac{Tr(G)}{2}}+1$.
		
		For the second case, if neither $x_{t+1}y_{t}\in E$ nor $x_{t}y_{t+1}\in E$, then $G$ contains $K_{t, t-1}$ as a subgraph for maximum $t$. Note that $x_ty_{t+i}\notin E$ and $x_{t+j}y_t\notin E$, for all $i, j\geq 2$. So, in this case, $G$ contains $(K_{2\floor{\frac{t'}{2}}, 2\floor{\frac{t'}{2}}-1}\setminus C_4)$ as a subgraph for maximum $t'$, where $t'=t$ if $t$ is odd and $t'=t+1$ if $t$ is even. Now, by the Lemma \ref{two_transitivity_chain_lemma_2}, we have $Tr_2(G)=\floor{\frac{t'}{2}}+1=\floor{\frac{t}{2}}+1=\floor{\frac{Tr(G)-1}{2}}+1$, when $t$ is odd. Otherwise, $Tr_2(G)=\floor{\frac{t'}{2}}+1=\floor{\frac{t+1}{2}}+1=\frac{t}{2}+1=\frac{Tr(G)-1}{2}+1=\floor{\frac{Tr(G)-1}{2}}+1$.

		Finally, for the third case, by the results of \cite{paul2023transitivity}, we have $Tr(G)=t+1$ and $K_{t, t}\setminus \{e\}$ is induced by $\{x_1, x_2, \ldots, x_t\}\cup \{y_1, y_2, \ldots, y_t\}$, where $e=x_ty_t$. Since $G$ contains $K_{t,t}\setminus \{e\}$ as an induced subgraph, $G$ contains $(K_{2\floor{\frac{t'}{2}}, 2\floor{\frac{t'}{2}}-1}\setminus C_4)$ as a subgraph for maximum $t'$, where $t'=t$ if $t$ is odd and $t'=t+1$ if $t$ is even. This follows from the chain ordering of $G$. Now, by the Lemma \ref{two_transitivity_chain_lemma_2}, we have $Tr_2(G)=\floor{\frac{t'}{2}}+1=\floor{\frac{t}{2}}+1=\floor{\frac{Tr(G)-1}{2}}+1$, when $t$ is odd. Otherwise, $Tr_2(G)=\floor{\frac{t'}{2}}+1=\floor{\frac{t+1}{2}}+1=\frac{t}{2}+1=\frac{Tr(G)-1}{2}+1=\floor{\frac{Tr(G)-1}{2}}+1$.
	\end{proof}
	
	In \cite{paul2023transitivity}, authors show that transitivity can be solved in linear time for a bipartite chain graph. Therefore, we have the following corollary of the above theorem.

	\begin{corollary}
		The \textsc{M$2$TP} can be solved in linear time for bipartite chain graphs.
	\end{corollary}

	\section{Conclusion}
	
	In this paper, we have introduced the notion of $2$-transitivity in graphs, which is a variation of transitivity. First, we have shown some basic properties for $2$-transitivity. We have shown that the \textsc{Maximum $2$-Transitivity Decision Problem} is NP-complete for chordal and bipartite graphs. On the positive side, we have proved that this problem can be solved in linear time for trees, split graphs, and bipartite chain graphs. It would be interesting to investigate the complexity status of this problem in other graph classes. Designing an approximation algorithm for this problem would be another challenging open problem.

	\section*{Acknowledgements:} Subhabrata Paul was supported by the SERB MATRICS Research Grant (No. MTR/2019/000528). The work of Kamal Santra is supported by the Department of Science and Technology (DST) (INSPIRE Fellowship, Ref No: DST/INSPIRE/ 03/2016/000291), Govt. of India.\\

	\bibliographystyle{alpha}
	\bibliography{Two_Transitivity_bibliography}

\end{document}